\documentclass[12pt]{article}
\usepackage[english]{babel}

\usepackage{graphicx}
\title{Curves on Rational Surfaces with Hyperelliptic Hyperplane Sections}
\author{Ovidiu P\u as\u arescu\footnote{Work partially supported by JSPS (Japan Soc. for Promotion of Science), Romanian Ministry of Education and Research (CERES Program) and Romanian Academy (GAR Project No. 5 / 5.02.2008)}}
\date{}
\begin{document}
\maketitle

\textheight 20cm
\textwidth 13cm
\oddsidemargin 12mm
\evensidemargin 4mm
\topmargin 10mm

\newcommand{\omm}{\Omega}
\newcommand{\om}{\omega}
\newcommand{\br}{\hbox{\bf R}}
\newcommand{\bn}{\hbox{\bf N}}
\newcommand{\bz}{\hbox{\bf Z}}
\newcommand{\bq}{\hbox{\bf Q}}
\newcommand{\bc}{\hbox{\bf C}}
\newcommand{\bp}{\hbox{\bf P}}
\newcommand{\bg}{\hbox{\bf G}}
\newcommand{\ba}{\hbox{\bf A}}
\newcommand{\pp}{\prime}
\newcommand{\ty}{\infty}
\newcommand{\di}{\displaystyle}
\newcommand{\va}{\varphi}
\newcommand{\si}{\sigma}
\newcommand{\ga}{\gamma}
\newcommand{\gaa}{\Gamma}
\newcommand{\na}{\nabla}
\newcommand{\te}{\theta}
\newcommand{\ld}{\ldots}
\newcommand{\ov}{\over}
\newcommand{\ri}{\Rightarrow}
\newcommand{\rii}{\Leftrightarrow}
\newcommand{\noa}{\noalign{\medskip}}
\newcommand{\la}{\lambda}
\newcommand{\su}{\subset}
\newcommand{\qu}{\quad}
\newcommand{\fo}{\forall}
\newcommand{\al}{\alpha}
\newcommand{\be}{\beta}
\newcommand{\ep}{\varepsilon}
\newcommand{\pa}{\partial}
\newcommand{\ti}{\times}
\newcommand{\dd}{\Delta}
\newcommand{\de}{\delta}
\newcommand{\sgg}{\Sigma}
\newcommand{\mm}{\medskip}
\newcommand{\hh}{\hfill q.e.d.}
\newcommand{\laa}{\langle}
\newcommand{\raa}{\rangle}
\newcommand{\ded}{\hbox{deg}}
\newcommand{\pic}{\hbox{Pic}}
\newcommand{\im}{\hbox{Im}}

\begin{abstract}
In this article we study, given a pair of integers $(d,g)$, the problem of existence of a smooth, irreducible, non-degenerate curve in $\bp^n$ of degree $d$ and genus $g$ (the Halphen-Castelnuovo problem). We define two domains from the $(d,g)$-plane, $D^n_1$ and $D^n_2$, and we prove that there is no gap in $D^n_1$. This follows constructing curves on some  rational surfaces with hyperelliptic hyperplane sections and from some previous theorems of Ciliberto, Sernesi and of the author. Moreover, in the last section, based on some results of Horrowitz, Ciliberto, Harris, Eisenbud, we conjecture that $D^n_2$ is the right lacunary domain.

\end{abstract}

\bigskip

\centerline{\Large \bf Contents}

\bigskip

1. Introduction \dotfill 2

2. The functions $\al_p (d,n)$; the domains $D^n_1, D^n_2$ and
$A^n_p$; the Main Theorem: absence of gaps in $D^n_1$ \dotfill 5

3. The proof of Main Theorem \dotfill 8

\hspace{1 cm} 3.1 Methods \dotfill 8

\hspace{1 cm} 3.2 Two smoothing criteria; the surfaces $X^n_p$ \dotfill 12

\hspace{1 cm} 3.3 Numerical properties of the functions $\al_p (d,n)$ \dotfill 17

\hspace{1 cm} 3.4 Invertible sheaves ${\cal L} \in \pic (X^n_p)$
whose $(\deg ({\cal L}), p_a ({\cal L}))_{\cal L}$ cover the domain
$A^n_p$ $(n\ge 5, p\ge n/3)$ \dotfill 18

\hspace{1 cm} 3.5 Curves $C\su X^n_p$ whose $(\deg (C), g(C))_C$
cover the domain $A^n_p$ $(n\ge 8, n/3 \le p \le n-4)$ \dotfill 28

\hspace{1 cm} 3.6 Conclusion \dotfill 36

 4. Comments and further developments; expected
gaps in $D^n_2$ \dotfill 36

References \dotfill 37

\section{Introduction}

In this article we work over an algebraically closed field $K$ of characteristic zero  (e.g. $K = \bc$ the complex field); see the section 4 also for arbitrary characteristic. We'll use the standard notations, see [Ha1].

By a {\it curve} (resp. a {\it surface}) we mean a $K$-algebraic integral scheme of dimension 1 (resp. 2). If not otherwise specified the curves and surfaces appearing in this paper will be supposed non-singular.

A (possible singular) curve $C\su \bp^n$ (= the projective $n$-space) is called {\it non-degenerate} if it is not contained in any hyperplane.

Given a property ${\cal P}$ we call ${\cal P}$-{\it Halphen Castelnuovo Theory in} $\bp^n$ the study of existence - for arbitrary, but fixed triplets if integers $(n,d,g)$, $n\ge 2$, $d\ge n$, $g\ge 0$ - of non-degenerate (possible singular) curves $C\su \bp^n$ of degree $d$ and geometric genus $g$, and having the property ${\cal P}$. Then it is natural to study the (non-empty) families ${\cal F}^{n,{\cal P}}_{d,g}$ of curves as before.

The property ${\cal P}$ could be, for instance, smoothness,
linear normality, projective normality, maximal rank, expected
number of moduli, etc. If ${\cal P}$=smootheness and $H^n_{d,g}$
is the Hilbert scheme of smooth curves from $\bp^n$ (i.e. the
closure in the general Hilbert scheme of the locus of smooth,
irreducible, non-degenerate curves from $\bp^n$), then one could
take ${\cal F}^{n{\cal P}}_{d,g} = H^n_{d,g}$.

When studying the non-degenerate {\it embedded curves} $C\su \bp^n$ for some $n\ge 2$ one usually study the previous schemes $H^n_{d,g}$, the first step being to find the triplets of integers $(n,d,g)$ for which $H^n_{d,g} \ne \emptyset$. Let's now denote, as usual, by ${\cal M}_g$ the moduli space of curves of genus $g$. Then, there is a natural rational map $\pi = \pi^n_{d,g}: H^n_{d,g} - - \to {\cal M}_g$ whose fibers over the $[C] \in {\cal M}_g$ are related to the determinantal varieties $W^n_d = W^n_d (C):= \{g^s_d| g^s_d$ a base-point free linear system of degree $d$ and dimension $s,s \ge n\}$ appearing in the Brill-Noether Theory for {\it polarized curves} $(C, g^s_d)$. As explained by J. Harris in [H], the rational map $\pi$ is very useful when studying the connections between the {\it extrinsic geometry} (represented  by the properties of projective embeddings and of the Hilbert scheme) and the {\it intrinsic geometry} (represented by the properties of the {\it abstract curves} and of the moduli space) of the families of curves. Certainly, this is possible when $\pi \ne \emptyset$ (so, again when $H^n_{d,g} \ne \emptyset$).

We recall now the {\it Castelnuovo bound} (found in 1893, for a modern proof, see [H], ch.3):

{\bf Theorem 1.1 (Castelnuovo [C])}. {\it If $n\in \bz$, $n\ge 3$ and $C\su \bp^n$ is a non-degenerate (possibly singular) curve of degree $d$ and geometric genus $g$, then $d\ge n$ and $0 \le g \le \pi_0 (d,n)$}.

This result has been a generalization of a previous similar result, obtained by Halphen and Noether in 1882 for $n=3$ ([Hl], [N]).

Here $\pi_0 (d,n)$ is the first one of the Harris-Eisenbud numbers $\pi_p (d,n)$, introduced in [H], ch.3 (for $p\ge 1$), and given (for $0\le p\le n-2$) by:
$$
\pi_p = \pi_p (d,n) = {m_p (m_p -1) \ov 2} (n+p-1) + m_p (\ep_p +p) + \mu_p , \qu \hbox{where} \leqno (1.1)
$$
$$
m_p = m_p (d,n) = [(d-1)/(n+p-1)]_* \leqno (1.2)
$$
(we denote, during this article, by $[x]_*$ the integer part of $x \in \br =$ the field of the real numbers)
$$
\ep_p = \ep_p (d,n) = d-1-m_p (n+p-1) \leqno (1.3)
$$
$$
\mu_p = \mu_p (d,n) = \max (0, [(p-n+2 + \ep_p)/2]_*). \leqno (1.4)
$$

We remark that $\mu_0 =0$ and $\pi_p = d^2/(2(n+p-1)) + O(d)$.

We recall, as an example, that if ${\cal P}$=nodality, the complete answer in ${\cal P}$-Halphen-Castelnuovo Theory is given by the following:

{\bf Theorem 1.2 (Tannenbaum $[T_1], [T_2]$)}. {\it For $n\ge 2$, $n\in \bz$ and any $d,g \in \bz$, $d\ge n$ and $0\le g \le \pi_0 (d,n)$ there is a non-degenerate curve $C\su \bp^n$ of degree $d$ and geometric genus $g$ having only nodes as singularities}.

This Theorem generalizes a similar result of Severi from 1915 for $n=2$ ([Sv]).

If ${\cal P}$ = smoothness we arrive to the {\it Halphen-Castelnuovo Problem} (related to the classification of the abstract curves and to the comparison between the intrinsic and extrinsic geometry of curves), denoted by $HC(n)$, for $\ge 2$, $n\in \bz$:

{\bf HC(n)}: For which pairs of integers $(d,g)$, $d\ge n$, $0\le g \le \pi_0 (d,n)$ do we have $H^n_{d,g} \ne \emptyset$ (equivalently, for which pairs $(d,g)$ as before there is a (smooth, irreducible) non-degenerate curve $C \su \bp^n$ of degree $d$ and genus $g$?).

This is the Problem which we'll consider in this article.

{\bf Definition 1.3.} {\it A pair of integers $(d,g)$, $d\ge n$, $0\le g \le \pi_0 (d,n)$ is called a {\bf gap} for $HC(n)$ if $H^n_{d,g} = \emptyset$ (i.e. there is no non-degenerate (smooth, irreducible) curve $C \su \bp^n$ of degree $d$ and genus $g$).}

From the Theorem 1.2 follows that there is no gap in ${\cal P}$-Halphen-Castelnuovo Theory when  ${\cal P}$ = nodality. However, the situation is different for $HC(n)$, so the Problem is more complicated.

Indeed, $HC(2)$ is simple: $H^2_{d,g} \ne \emptyset \rii d\ge 2$ and $g = (d-1)(d-2)/2$. But, for $n\ge 3$, $HC(n)$ becomes highly non-trivial. A (correct) answer to $HC(3)$ has been stated by  Halphen in 1882 ([Hl]), but his proof has been, partially, incorrect. A complete proof has been given, 100 years later, by Gruson and Peskine ([GP1], [GP2]). The solution of $HC(3)$ is the following:

{\bf Theorem 1.4 (Halphen-Gruson-Peskine)}. {\it For $HC(3)$ there are two domains in the $(d,g)$-plane: the {\bf non-lacunary} domain $D^3_1$, where there is no gap for $HC(3)$ and the {\bf lacunary domain} $D^3_2$, where there are gaps. Here
$$
\begin{array}{l}
D^3_1: 0\le g \le \pi_1 (d,3), \; d\ge 3 \\ \noa
D^3_2: \pi_1 (d,3) < g \le \pi_0 (d,3), \; d\ge 3. \end{array}.
$$
Moreover, a pair $(d,g) \in D^3_2$ is {\bf not} a gap for $HC(3)$ iff $(d-2)^2-4g$ is the square of some integer (see the Theorem 1.1 also)}.

$HC(n)$ has been also solved for $n=4,5$ by Rathmann ([Ra]), $n=6$ by Ciliberto ([Ci]) and "almost" solved by Ciliberto if $n=7$ ([Ci]). The situation is similar to the case $n=3$ in the sense that there are two domains: $D^n_1$, where there is no gap, and $D^n_2$, where there are gaps. The domains $D^n_1$ and $D^n_2$ $(n \in \{4,5,6,7\})$ are explicitely defined and the gaps from $D^n_2$ are also determined (with the exception of 24 pairs for $n=7$ for which is not yet known if they are gaps or not). Results for the general $HC(n)$ belongs to Ciliberto, Sernesi ([Ci], [CS]), Harris, Eisenbud ([H]), Rathmann ([R]), P\u as\u arescu ([P1], [P2], [P3]) and others. We'll use some of them in the next section (see Section 4 also).

We end this section remarking that  the case $n=3$ of $HC(n)$ has some intrinsic importance, because any (smooth, irreducible) curve from some $\bp^n$, $n\ge 4$, can be projected (so, preserving the degree) isomorphically (so, preserving the genus) on $\bp^3$. But, for more sophisticated reasons, consisting in the comparison between the extrinsic and intrinsic properties, the natural problem is the general $HC(n)$, $n\ge 3$.

{\bf Aknowledgement}: During the preparation of this article, I have worked to the Institute of Mathematics of the Romanian Academy and I have visited the University of Ferrara, the Mittag-Lehfler Institute from Stockholm, the University of Nice, the University of Angers, TIFR Mumbai, MSRI Berkeley, the National Australian University Canberra, the University of Iowa and the University of Tokyo (in this chronological order). I have had interesting and stimulating discussions on the subject with J. Alexander, C. Ciliberto, J. Coand\u a, D. Eisenbud, Ph. Ellia, R. Hartshorne, A. Hirschowitz, T. Katsura, J. Kleppe, D. Laksov, A. Lascu, Y. Miyaoka, C. Peskine, M. Reid, R. Strano. I want to thank to \c S. Basarab, A. Buium, A. Dimca, D. Popescu also for general stimulating discussions.

\section{The functions $\al_p (d,n)$; the domains $D^n_1, D^n_2$ and $A^n_p$; the Main Theorem: absence of gaps in $D^n_1$}

Before stating the Main Theorem which will be proved in this article, we need a number of definitions.

For $p,n \in \bz$, $p\ge 0$, $n\ge 3$ we {\it define} the following numerical functions $\al_p (d,n) = d^2/(2(n+p-1)) + O(d)$ given by:
$$
\al_p = \al_p (d,n) = {x_p (x_p-1) \ov 2} (n+p-1) + x_p (t_p+p) + u_p,
\qu \hbox{where} \leqno (2.1)
$$
$$
x_p = x_p (d,n): = [(d-a^n_p)/ (n+p-1)]_*  \leqno (2.2)
$$
$$
a_p = a^n_p := [(n-p)/2]_* + 1  \leqno (2.3)
$$
$$
t_p = t_p (d,n): = d-1-x_p (n+p-1)  \leqno (2.4)
$$
$$
u_p = u_p (d,n): = [(p-n+1+t_p)/2]_* \; .  \leqno (2.5)
$$

We also consider the functions
$$
d_1(n):= \left\{\begin{array}{ll}
\max (2n+1, \di{1\ov 6} (3 + (4k-1) \sqrt{24k - 33})), & n\equiv 0 \; (\hbox{mod 3}) \\ \noa
\max (2n+1, \di{4k\ov 4k+1} (-4k+ \sqrt{32k^3 + 16k^2 -2k-1})), & n\equiv 1 \; (\hbox{mod 3}) \\ \noa
\max (2n+1, 5k+ 3 + (2k + 1) \sqrt{48k + 6}), & n\equiv 2 \; (\hbox{mod 3}) \end{array} \right. \leqno (2.6)
$$
where $k= [n/3]_*$, $n\ge 4$.

We'll need also (with $k$ from before):
$$
A(d,n):= \left\{\begin{array}{ll}
\pi_k (d,n), & \hbox{if} \; n \equiv 0,1 \; (\hbox{mod 3}) \\ \noa
\pi_k (d,n) - \mu_k + \max (0, \ep_k -3k-1), & \hbox{if} \; n\equiv 2 \; (\hbox{mod 3}) \end{array} \right. \leqno (2.7)
$$
(see (1.1)-(1.4)).

\medskip

$$
B(d,n):= \left\{\begin{array}{ll} \left. \begin{array}{ll} \al_{k+1}
(d,n), & \hbox{if} \; n\equiv 1,2 \; (\hbox{mod 3}) \\ \noa \al_k
(d,n), & \hbox{if} \; n\equiv 0 \; \hbox{(mod
3)}\end{array}\right\}, & \hbox{if} \; 2n+1 \le d < d_1 (n)\\ \noa
A(d,n), & \hbox{if} \; d \ge d_1 (n) \end{array} \right. \leqno
(2.8)
$$

We consider now the following graphs (contained in the $(d,g)$-plane), of equations:
$$
\begin{array}{l}
B^n_0: \; g = \pi_0 (d,n) \\ \noa
B^n_k: \; g = \left\{\begin{array}{ll} \pi_0 (d,n), & n\le d\le 2n \\ \noa
B(d,n), & d\ge 2n+1 \end{array}\right., \; k = [n/3]_*  \\ \noa
C^n_p: \; g = \al_p (d,n), \; p\ge n/3, \; p,n \in \bz, \; n\ge 3 \end{array}
$$

Now we are ready to {\it define} the domains $D^n_1$ (bounded by $g=0$ and $B^n_k$) and $D^n_2$ (bounded by $B^n_k$ and $B^n_0$), contained in the $(d,g)$-plane:
$$
D^n_1: 0\le g \le \left\{\begin{array}{ll} \pi_0 (d,n), & n\le d\le 2n \\ \noa
B(d,n), & d\ge 2n+1 \end{array}\right. \leqno (2.9)
$$
$$
D^n_2: B(d,n) < g \le \pi_0 (d,n), \; d\ge 2n+1 \leqno (2.10)
$$

The Main Theorem belonging to Halphen-Castelnuovo Theory which we'll prove in this article is:

{\bf MAIN THEOREM}. {\it In the domain $D^n_1$ there is no gap for $HC(n)$, $n\ge 4$, $n\in \bz$.}

{\bf Remark 2.1.} {\it The Main Theorem is known to be true for $d> D(n) =$ a quadratic function in $n$ (see} [Ci], [CS], [P1]). {\it So, our main contribution is for "small" degrees, $2n+1 \le d < D(n)$. Moreover, as we recalled in section 1 the Main Theorem  is known to be true for $n \in \{4,5,6,7\}$, due to the contributions of Rathmann} ([Ra]), {\it Ciliberto, Sernesi} ([Ci], [CS]), {\it and of the author} ([P1]). {\it So, we'll suppose that $n\ge 8$ (but our arguments work, in principle, for $5 \le n \le 7$ also).}

During the proof of the Main Theorem, the domain $D^n_1$ will be divided in subdomains $A^n_p$, defined by:
$$
A^{'n}_{n-3}: 0 \le g \le \al_{n-3} (d,n), \; d\ge 2n+1 \; (\hbox{between $g=0$ and $C^n_{n-3}$)} \leqno (2.11)
$$
$$
A^n_p: \al_{p+1} (d-1,n) \le g \le \al_p (d,n), \; d\ge 2n+1, p \ge n/3 \leqno (2.12)
$$
(containing the domain between $C^n_{p+1}$ and $C^n_p$ for $p\le n-1$: see subsection 3.3 from section 3).
$$
\tilde A^n_k: \al_{k+1} (d-1,n) < g \le B (d,n), \; d\ge 2n+1, k = [n/3]_* \leqno (2.13)
$$
(containing the domain between $C^n_{k+1}$ and $B^n_k$)
$$
A^n_0: 0\le g \le \pi_0 (d,n), \; n\le d \le 2n. \leqno (2.14)
$$

{\bf Lemma 2.2.} {\it In the domain $A^n_0$ there is no gap for $HC(n)$, $n\ge 8$}.

{\bf Proof}. We use projections of non-special curves. For $d=2n$ and $g = \pi_0(d,n) = \pi_0 (2n,n)$, we consider the extremal curve.

\hh

The following equality holds (see Lemma 3.3.1b)):
$$
D^n_1 = A^n_0 \cup A^{'n}_{n-3} \cup \left( \bigcup_{{n\ov 3} \le p \le n-4} A^n_p\right) \cup \tilde A^n_k \leqno (2.15)
$$

We just solved the elementary step (curves in $A^n_0$) of $HC(n)$, $n\ge 8$ in the previous lemma. By (2.15), the Main Theorem is now the consequence of the following three Theorems:

{\bf Theorem A}. {\it If $n\in \bz$, $n\ge 8$, then there is no gap for $HC(n)$ in the domain $\tilde A^n_k$}.

The {\it proof} of this Theorem essentially belongs to Ciliberto ([Ci], Theorem 1.1). However the statement of the Theorem A differs slightly from the statement of the Theorem of Ciliberto. The proof of the Theorem A can be found in [P2].

{\bf Theorem B}. {\it If $n\in \bz$, $n\ge 8$, then there is no gap for $HC(n)$ in the domain $A^{'n}_{n-3}$}.

The {\it proof} of this Theorem essentially belongs to Ciliberto and Sernesi ([CS], Main Theorem), see also [P1]. However the statement of the Theorem of Ciliberto and Sernesi differs slightly from the statement of the Theorem B. The proof of the Theorem B can be found in [P3].

{\bf Theorem C}. {\it If $n\in \bz$, $n\ge 8$, then there is no gap for $HC(n)$ in the domain $\di\bigcup_{{n\ov 3} \le p \le n-4} A^n_p$}.

This Theorem C represents our main contribution. It has been already used during the proof of the Theorem A in [P2]. The complete proof of the Theorem C is given in the next section 3. It is done constructing curves on some rational surfaces $X^n_p \su \bp^n$ with hyperelliptic hyperplane sections.

%%%%We end this section with a picture:

%%%$$\begin{array}{c}
%%%%{\includegraphics*[width=12cm,height=7cm]{ovidiu-2.eps}}\\
%%%%\hbox{Fig. 1. }\\
%%%%\end{array}$$

%%%The hashed domain in $D^n_1$ (a non-lacunary domain).

\section{The proof of Main Theorem}

\subsection{Methods}

Let's fix some notations. Let there be $\Sigma \su \bp^2$ a finite set of (distinct) points $\Sigma = \{P_0, P_1, \ld , P_s\}$. We denote by $S = BL_{\Sigma} (\bp^2) \to \bp^2$ the blow up of $\bp^2$ in $\Sigma$. Then $\pic (S) \cong \bz \oplus \bz^{s+1}$ with $(l; -e_0, -e_1, \ld, -e_s)$ a $\bz$-basis (here $l$ is the class of the inverse image in $S$ of a line $L \su \bp^2$ and $e_i$ are the classes of the exceptional divisors $E_i \su S$ corresponding to the points $P_i$, $0\le i\le s$). We recall that the intersection form on $\pic (S)$ is given by:
$$
(l\cdot e_i) = 0, \; (e^2_i)=-1, \; i = \overline{0,s}, \; (e_i\cdot e_j) =0, \; (\fo) i\ne j, \; (l^2)=1.
$$

If ${\cal L} = al - \di\sum^s_{i=0} b_i e_i \in \pic (S)$, we write ${\cal L} = (a; b_0,b_1, \ld, b_s)$.

If ${\cal L} = {\cal O}_S (D)$, $D\in \hbox{Div}\: (S)$, we denote by $[{\cal L}] = |D|$ the complete linear system associated to ${\cal L}$ (resp. $D$). We write
$[(a; b_0,b_1, \ld, b_s)]:= [a; b_0,b_1, \ld, b_s]$. We say that $a,b_0,b_1, \ld, b_s$ are the {\it coefficients} of ${\cal L}$ (or [${\cal L}$]).

Let's now analyse the Gruson-Peskine type construction given by Ciliberto ([Ci]). The surfaces $X^n_k$ used by Ciliberto, containing the necessary curves $(k = [n/3]_*)$ are obtained blowing up $3k-n+6$ points from $\bp^2$ in {\it general position} (i.e. any 3 non-collinear and not all of them on a smooth conic) and embedding the abstract surface obtained in such a way using the very ample invertible sheaf $(k+2;k,1^{3k-n+5})$ (we denoted $\di(a;b, \underbrace{c, \ld,c)}_{s\: times}$ by $(a; b, c^s$)). Because $3k-n+6 \in \{4,5,6\}$ there is an well-known criterion giving sufficient conditions for a linear system in order to contain (smooth, irreducible) curves. Precisely, if
${\cal L} = (a; b_0,b_1, \ld, b_{s^n_k}) \in \pic (X^n_k)$, $s^n_k = 3k-n+5$ and
$$
a\ge b_0 + b_1 + b_2, \; b_0 \ge b_1 \ge \ld \ge b_{s^n_k} \ge 0, \; a> b_0 \leqno (3.1.1)
$$
then $[{\cal L}] = [a; b_0,b_1, \ld, b_{s^n_k}]$ contains a (smooth, irreducible) curve.

Due to this criterion (depending only on coefficients), the construction of curves on the surfaces $X^n_k$ has two parts: an {\it arithmetical part}, where the necessary degrees and arithmetical genera are obtained using the well-known formulae of genus and degree for various ${\cal L} \in \pic (X^n_k)$, and a {\it geometrical part}, where the necessary (smooth, irreducible) curves are obtained applying (3.1.1) to the sheaves ${\cal L}$ used in order to solve the arithmetical part.

The main contribution of this paper to $HC(n)$ is the construction of invertible sheaves (and then curves) having the degree and genus in the domains $A^n_p$ (see (2.12)), on some surfaces $X^n_p$. Precisely, let's denote by $S_p = S^n_p: = Bl_{\Sigma_p} (\bp^2)$,
$\Sigma_p = \Sigma^n_p = \{P_0, P_1, \ld, P_{s^n_p}\} \su \bp^2$, $s^n_p = 3p -n+5$,
$\Sigma^n_p$ containing general points. Then ${\cal H}^n_p:= (p+2; p, 1^{s^n_p}) \in \pic (S^n_p)$ is very ample (see subsection 3.2, Prop. 3.2.2) and then $X^n_p:= \im \va_{[{\cal H}^n_p]} \su \bp^n$, $\deg X^n_p = n+p-1$. If $p=k$ one obtains the surfaces $X^n_k$ used in [Ci].

If we try to use an argument similar to the argument from  [GP2] or [Ci] for $X^n_p$, so to divide the proof in two parts, an arithmetical and a geometrical one (the  smoothing) we need a criterion similar to (3.1.1) for an {\it arbitrary} number of points from $\Sigma_p$. But this is, anyway, very complicated by itself (see, for instance [Hi]). On the other hand, such a general criterion must include the conditions (3.1.1) ({\it normalization}, using quadratic transformations, see [Hi]); but these conditions allow us to construct only curves whose degree $d$ is $d>a$ function of degree 3/2 in $n$, as it can be seen if we try to apply an argument as in [Ci]. So, in order to produce the necessary curves {\it for small degrees}, we need a new technique, which we shall explain here.

The main idea is to consider the arithmetical and the geometrical parts of the construction of curves on $X^n_p$ {\it entirely linked}, i.e. to construct linear systems containing the necessary curves directly and explicitly enough using combinations between some simple sheaves (so, somehow easy to understand them).

Explicitly, we'll proceed as follows (in order to construct the necessary curves on $X^n_p$ whose (degree, genus) belongs to the domain $A^n_p$):  we start with a simple {\it initial family} ${\cal D}_0$ of sheaves
$$
{\cal D}_0 = (a+2; a, 1^t, 0^{s^n_p-t}) \in \pic (X^n_p), \; 0\le t\le s^n_p, \; a \in \bz \leqno (3.1.2)
$$
of arithmetical genus $a$ and some degrees,  realizing the necessary genera in some initial intervals in $d$ (like $x_p (d,n)=0$, for instance - see (2.2) - but not only). The families (3.1.2) are the only one which are good for the initial intervals in our method (we'll explain this later). After the initial construction we continue using a number of inductive arguments (after $x_p (d,n)$, for instance, but not only), adding repeatedly to ${\cal D}_0$ (by tensorization) some simple (and well understood) invertible sheaves, similar to the class of hyperplane section (see Lemma 3.2.4). Let's explain, shortly, why the inductive processes works. Let's suppose that we need to construct curves on the surfaces $X^n_p$ of degree $d$ and genus $g = \pi_p (d,n)$ (see (1.1)). Let's suppose that we already constructed the necessary curves for $m_p (d,n)=0$ (see (1.2)). Because
$$
\begin{array}{c}
\pi_p (d+(n_p-1)) = \pi_p (d) + (d+p-1); \\ \noa
 m_p (d+(n+p-1)) = m_p (d)+1  \end{array} \leqno (3.1.3)
$$
in order to construct curves $C$ with $(\deg (C), g(C)) = (d,g)$,
we proceed as follows: if the curves $C$ have degrees $d$ and
genera $g = \pi_p (d,n)$ in the domain $m_p (d,n)=m$, then the
(smooth, irreducible!) curves $C' \in |C + H^n_p|$, $H^n_p \in
[p+2; p, 1^{s^n_p}]$ have degrees $d' = d+(n+p-1)$ and genera $g'
= \pi_p (d',n)$ in the next domain $m_p (d,n) =m+1$ (use (3.1.3)).
Hence, we obtain invertible sheaves for all pairs $(d,g)$, $g =
\pi_p (d,n)$ as far as we succeed to construct these invertible
sheaves in the domain  $m_p (d,n)=0$. Moreover, we can test when
the associated linear system contain (smooth, irreducible) curves,
because the appearing linear systems $|D_0 + t H^n_p|$, $t \ge 0$
are simple enough (specialize the points from $\Sigma^n_p$ on a
rational plane curve of degree $p+1$, having one ordinary singular
point $P_0$ with multiplicity $p$,  see Prop. 3.2.2 and Lemma
3.2.4).

In our case, in order to construct the necessary curves covering by (degree, genus) the domains $A^n_p$ it is necessary to change in some way the functions $\pi_p$ so that the initial verifications (corresponding to $m_p (d,n)=0$) to be made using sheaves as in (3.1.2) and a property as in (3.1.3) to still hold (and some others appearing, see subsection 3.3, Lemma 3.3.1). We get in such a way the (unique) functions $\al_p (d,n)$ from section 2 (see (2.1)).

Moreover, the method shortly explained here works only if $\# \Sigma^n_p = 3p-n+6\ge 12$ ($\# A$ means the number of elements of the finite set $A$) so, if $p\ge n/3+2$. So, for $\# \Sigma^n_p \le 11$ we need {\it another method} in order to construct curves on $X^n_p$. Because we need such curves for {\it small degrees} also, the inductive argument is used again, as possible. But in this second method of construction we'll be, partially inspired from [GP2], because in this case a smoothing criterion of type (3.1.1) will be good enough (it is obtained specializing the points from $\Sigma^n_p$ on a smooth plane cubic curve, see Prop. 3.2.1 and 3.2.3). In this case, a {\it supplementary property} of the functions $\al_p$ will be necessary, namely that these functions must be related in a "good" way with the genus formula (see subsection 3.3, Lemma 3.3.2); and this is possible exactly because one uses initial families of type (3.1.2) (so these families are obligatory).

The smoothing criteria used in the methods are collected in the subsection 3.2, in Prop. 3.2.3 and Lemma 3.2.4. The numerical  properties of the functions $\al_p$ are presented in the subsection 3.3. In the subsection 3.4 we'll construct invertible sheaves from $\pic (X^n_p)$ covering by (degree, arithmetical genus) the domain $A^n_p$, using both methods. In subsection 3.5 we'll apply the two smoothing criteria to sheaves from the subsection 3.4 in order to get (smooth, irreducible) curves.  We remark that the domain of applicability of the first method  briefly described here is (R1), while the domain of applicability of the second method (using some of the Gruson-Peskine ideas) is (R2), where:
$$
(R1): n \ge 9, \; n/3 + 2 \le n-4 \; (\hbox{so}\; \# \Sigma^n_p \ge 12) \leqno (3.1.4)
$$
$$
(R2): n \ge 8, \; n/3 \le p , n/3 + 2 \; (\hbox{so}\; \# \Sigma^n_p \le 11) \leqno (3.1.5)
$$

We end this section remarking that the constructions from the subsection 3.4 and 3.5 have {\it no degree of freedom}, but they cover all the ranges. Let's remark that the  method used in (R2) cannot be used in general, because the condition (C5) from Prop. 3.2.3  gives   (smooth, irreducible) curves only for $d>a$ quadratic function in $n$ if used for $p=n-4$, $n-5, \ld$

\subsection{Two smoothing criteria; the surfaces $X^n_p$}

Let there be $\Sigma \su \bp^2$, $\Sigma = \{P_0, P_1, \ld, P_s\}$ a set of general points and $S:= Bl_{\Sigma} (\bp^2) \stackrel{\pi}{\to} \bp^2$ the blow up of $\bp^2$ in $\Sigma$.

The first smoothing criterion is ii) from the next Proposition 3.2.1 (reformulated in Prop. 3.2.3).

{\bf Proposition 3.2.1}. {\it If ${\cal D} = (a; b_0, b_1, \ld, b_s)\in \pic (S)$ is so that $a\ge b_0 \ge b_1 \ge \ld \ge b_s >0$, then:

i) each one of the conditions (i.1), (i.2), (i.2)$'$, (i.3), (i.3)$'$ implies $h^1 ({\cal D}) =0$, where

(i.1): $a - \di\sum^s_{l=0} b_l \ge -1$;

(i.2): $2\le s \le 7$, $a \ge \di\sum^2_{l=0} b_l$;

(i.2)$'$: $2\le s \le 7$, $a - \di\sum^s_{l=0} b_l \ge -1$, $b_0 > b_1$;

(i.3): $s\ge 8$, $a \ge \di\sum^2_{l=0} b_l$, $3a -
\di\sum^s_{l=0} b_l \ge 1$;

(i.3)$'$: $s\ge 8$, $a - \di\sum^2_{l=0} b_l \ge -1$, $3a - \di\sum^s_{l=0} b_l \ge 1$, $b_0 > b_1$;

ii) each one of the conditions (ii.1), (ii.2), (ii.2)$'$, (ii.3), (ii.3)$'$ implies $h^0 ({\cal D})\ne 0$, $[{\cal D}]$ has no base point and contains a (smooth, irreducible) curve, where

(ii.1): $a\ge \di\sum^s_{l=0} b_l$;

(ii.2): $2\le s\le 6$, $a \ge \di\sum^2_{l=0} b_l$;

(ii.2)$'$: $2 \le s \le 7$, $a - \di\sum^2_{l=0} b_l \ge -1$, $b_0 \ge b_1 +2$;

(ii.3): $s\ge 7$, $a\ge \di\sum^2_{l=0} b_l$, $3a - \di\sum^s_{l=0} b_l \ge 2$;

(ii.3)$'$: $s\ge 8$, $a - \di\sum^2_{l=0} b_l \ge -1$, $3a - \di\sum^s_{l=0} b_l \ge 2$, $b_0 \ge b_1 + 2$.}

{\bf Remark}. This is a Harbourne type result ([Hb1]).

{\bf Proof}. A direct proof can be obtained specializing the points from $\Sigma$ on a smooth plane cubic $\gaa_0 \su \bp^2$ (general on $\gaa_0$), using induction on $s$ and standard exact sequences. The details are left to the reader (or see [P4]).

\hh

{\bf Proposition 3.2.2}. {\it Let there be ${\cal D} = (p+2; p, 1^s) \in \pic (S)$, $p\in \bz$, $p\ge 1$. Then: i) if $s \le 3p+3$, then $h^\circ ({\cal D}) \ne 0$, $[{\cal D}]$ has no base point and contains a (smooth, irreducible) curve;

ii) if $s\le 3p$, then ${\cal D}$ is very ample.}

{\bf Remark}. A proof of ii) can be found in [Gi].

{\bf Proof.} A direct proof can be obtained specializing  the points from $\Sigma$
(except for the last one) on a rational irreducible plane curve  $\dd_0 \su \bp^2$ of
degree $p+1$ having only one ordinary singularity of multiplicity $p: P_0$ will be
the singular point, $P_1, \ld, P_{s-1} \in \dd_0$, $P_s \notin \dd_0$. Then use
standard exact sequences and standard techniques in order to separate points and tangent
vectors. The details are left to the reader (or see [P4]).

\hh

Now, we are ready to {\it define} the surfaces $X^n_p$. Precisely,
let there be $p\in \bz$, $k\le p \le n-4$, $k = [n/3]_*$,
$n\in\bz$, $n\ge 5$ and $\Sigma_p = \Sigma^n_p: = \{P_0, P_1, \ld,
P_{s^n_p}\} \su \bp^2$ $(s^n_p:= 3p-n+5$) be a set of general
points. From Prop. 3.2.2 ii) follows that  ${\cal H}^n_p:= (p+2;
p, 1^{s^n_p}) \in \pic (S^n_p)$ $(S^n_p:= Bl_{\Sigma^n_p}
(\bp^2))$ is very ample on $S^n_p$. Then:
$$
X^n_p:= \im \left(\va_{[{\cal H}^n_p]}\right).
$$

The surfaces $X^n_p$ are rational having hyperelliptic hyperplane sections and one can easily check (doing computations) that $X^n_p \su \bp^n$ and $\deg (X^n_p) = n+p-1$. If $p=k$ $(= [n/3]_*)$ one obtains the surface used in [Ci].

Now we need to reformulate Prop. 3.2.3 ii) in the Gruson-Peskine coordinates
$(d,r: \te_1, \ld, \te_{s^n_p})$ ([GP2, step 2]); this {\it is possible}. So, for ${\cal D}\in \pic (X^n_p) \cong \bz \oplus \bz^{s^n_p+1}$, ${\cal D} = {\cal O}_X (D) = (a; b_0, b_1, \ld, b_{s^n_p})$, satisfying $a\ge b_0 \ge b_1 \ge \ld \ge b_{s^n_p} \ge 0$, $(D \in \hbox{Div} (X^n_p))$, we consider the change of coordinates (in $\pic (X^n_p) \otimes \bq)$ given by:
$$
r:= a-b_0, \; \theta_i:= {1\ov 2} r - b_i, \; i = \overline{1, s^n_p}. \leqno (3.2.1)
$$

If $d = \deg (D):= ({\cal D} \cdot H^n_p)$ $(H^n_p \in [{\cal H}^n_p])$ and $g = p_a (D)$, using the genus formula one obtains:
$$
d = {p-n+5 \ov 2} r+2a + \sum^{s^n_p}_{i=1} \theta_i; \; g = F_d (r) - {1\ov 2}
\sum^{s^n_p}_{i=1} \theta_i^2. \leqno (3.2.2)
$$
$$
F_d (r) = F^{p,n}_d (r):= {1\ov 2} \left[d(r-1)+(p-1)r - {n+p-1 \ov 4}r^2\right] +1. \leqno (3.2.3)
$$

Prop. 3.2.1 ii) can be reformulated now as:

{\bf Proposition 3.2.3}. {\it Let there be}
$$
{\cal D} = {\cal O}_{X^n_p} (D) \in \pic (X^n_p), \qu {\cal D} = (a; b_0,b_1, \ld, b_{s^n_p}).
$$

{\it If the following conditions (C1), (C2), (C3), (C4), (C5) are
simultaneously satisfied by ${\cal D}$ in the coordinates $(d,r;
\theta_1 , \ld, \theta_{s^n_p})$, then $[{\cal D}] = |D| \ne
\emptyset$, has no base point and contains a (smooth, irreducible)
curve, where}:
$$
\theta_i \equiv {1\ov 2} r \qu \hbox{(mod 1)}, \qu i = \overline{1, s^n_p}; \leqno (C1)
$$
$$
d + {1\ov 2} (p-n+5)r - \sum^{s^n_p}_{i=1} \theta_i  \equiv 0 \qu \hbox{(mod 2)} ; \leqno (C2)
$$
$$
|\theta_1| \le  \theta_2 \le \ld \le \theta_{s^n_p} \le {r\ov 2}; \leqno (C3)
$$
$$
-\te_1 + \sum^{s^n_p}_{i=2} \theta_i  \le d - {1\ov 2} (n-p+1)r; \leqno (C4)
$$
$$
d \ge (p-1)r + 2.  \leqno (C5)
$$

{\bf Proof}. Left to the reader (similar to the case $n=3$, see [Ha2], [GP2]).

\hh

We end this section with a technical lemma (which represents the second smoothing criterion), essential in section 3.5.

{\bf Lemma 3.2.4.} {\it If $p,n \in \bz$, $n\ge 3$, $p\ge {n\ov 3}
+2$, let's consider the following invertible sheaves from} $\pic
(X^n_p)$:
$$
\begin{array}{l}
{\cal H}_1 = {\cal H}^n_p: = (p+2; p, 1^{s^n_p}), \\ \noa
{\cal H}_2 = \tilde {\cal H}^n_p: = (p+2; p, 1^{s^n_p-1},0), \\ \noa
{\cal H}'_3 = (\tilde{\cal H}^{n+1}_{p-1,1})': = (p+1; p-1, 1^{s^n_p-5},0^2,1^2,0), \\ \noa
{\cal H}'_4 = (\tilde{\cal H}^{n+1}_{p-1,2})': = (p+1; p-1, 1^{s^n_p-7},0^2,1^4,0), \\ \noa
{\cal H}'_5 = (\tilde{\cal H}^{n+1}_{p-1,3})': = (p+1; p-1, 1^{s^n_p-9},0^2,1^6,0), \\ \noa
{\cal H}'_6 = (\tilde{\cal H}^{n+1}_{p-1,4})': = (p+1; p-1, 1^{s^n_p-11},0^2,1^8,0). \end{array}
$$

{\it If $t_1, t_2, \ld, t_6 \in \bz$, $t_1, \ld, t_6 \ge 0$, let
there be}
$$
{\cal D}:= {\cal D}_0 + t_1 {\cal H}_1 + t_2 {\cal H}_2 + t_3 {\cal H}'_3 + t_4 {\cal H}'_4 + t_5 {\cal H}'_5 + t_6 {\cal H}'_6 \in \pic (X^n_p) \leqno (3.2.4)
$$
{\it where
$$
{\cal D}_0:=  (a+2; a, \nu_1, \nu_2, \ld, \nu_{s^n_p}), \; \nu_j \in \{0,1\}, \qu
j = \overline{1, s^n_p}
$$
and
$\nu_j=1$ for $u$ values of the index $j$ $(0 \le u \le s^n_p)$. Let's denote by
$t:= t_1 + t_2 + \ld + t_6$. Then, the condition
$$
3a \ge u + n-7 - (t+1)(n-4) \leqno (3.2.5)
$$
implies $[{\cal D}] \ne \emptyset$, $[{\cal D}]$ without base points and contains a (smooth, irreducible) curve.}

{\bf Proof}. Write
$$
{\cal D} = t_1 {\cal D}_{x_1} + t_2 {\cal D}_{x_2} + t_3 {\cal D}_{x_3} +
t_4 {\cal D}_{x_4} + t_5 {\cal D}_{x_5} + t_6 {\cal D}_{x_6} + {\cal R}
$$
$x_1, \ld, x_6 \in \bz$, where
$$
\begin{array}{l}
{\cal D}_{x_1}:= (x_1 + 2; x_1, 1^{s^n_p}), \\ \noa {\cal
D}_{x_2}:= (x_2 + 2; x_2, 1^{s^n_p-1},0), \\ \noa {\cal D}_{x_3}:=
(x_3 + 2; x_3, 1^{s^n_p-5},0^2, 1^2,0), \\ \noa {\cal D}_{x_4}:=
(x_4 + 2; x_4, 1^{s^n_p-7},0^2,1^4,0), \\ \noa {\cal D}_{x_5}:=
(x_5 + 2; x_5, 1^{s^n_p-9},0^2,1^6,0), \\ \noa {\cal D}_{x_6}:=
(x_6 + 2; x_6, 1^{s^n_p-11},0^2,1^8,0), \\ \noa
\end{array}
$$
${\cal R}:= (b+2; b, \nu_1, \nu_2, \ld, \nu_{s^n_p})$, where $b:= a + t(p-1) + (t_1 + t_2) - \di\sum^6_{i=1} t_i x_i$.

Take then $x_1, \ld, x_6 \ge 0$ {\it minimal} so that  Prop. 3.2.2 i) applies to ${\cal D}_{x_1}, \ld, {\cal D}_{x_6}$. We deduce that:
$$
3p-n+2 \le 3x_1 \le 3p-n +4 \leqno (3.2.6)
$$
$$
3p-n+1 \le 3x_2 \le 3p-n +3 \leqno (3.2.7)
$$
$$
3p-n-1 \le 3x_i \le 3p-n +1, \qu i=\overline{3,6} \leqno (3.2.8)
$$

In order to obtain  the conclusion from the lemma for $[{\cal D}]$, it's enough to have the same for $[{\cal R}]$ so, by Prop. 3.2.2 i) again, we need
$$
u\le 3b+3, \; b = a+ t(p-1)+(t_1 + t_2) - \sum^6_{i=1} t_ix_i \leqno (3.2.9)
$$

Replacing $b$, (3.2.9) becomes
$$
u\le 3a + 3t (p-1) + 3(t_1 + t_2) - \sum^6_{i=1} t_i (3x_i) + 3, \qu \hbox{or}
$$
$$
\sum^6_{i=1} t_i (3x_i) \le 3a + 3t(p-1)+3(t_1+t_2) -u+3.
$$
Using (3.2.6), (3.2.7), (3.2.8) we can see that this last inequality is implied by
$$
t(3p-n+1)+3t_1+2t_2 \le 3a + 3t (p-1) + 3(t_1+t_2) -u+3 \rii
$$
$$
\rii 3a \ge u + n -7 -(t+1)(n-4)-t_2.
$$
Because $t_2 \ge 0$, this last inequality follows from (3.2.5).

\hh

{\bf Corollary 3.2.5}. {\it In the same hypotheses as in Lemma 3.2.4, the same conclusion holds for $[{\cal D}]$ if }
$$
3a \ge 3p-(t+1)(n-4)-2. \leqno (3.2.10)
$$

{\bf Proof}. Give to $u$ the biggest possible value $u = s^n_p = 3p-n+5$ in (3.2.5).

\hh

\subsection{Numerical properties of the functions $\al_p (d,n)$}

We recall that the functions $\al_p = \al_p (d,n)$ were defined in section 2 ((2.1)-(2.5)). Let it be now
$$
\al'_p (d,n):= \al_{p-1} (d+1, n+1). \leqno (3.3.1)
$$

We prove now the following {\it key lemma}:

{\bf Lemma 3.3.1}. {\it a) Let's denote by $d': =d+(n+p-1)$. Then:
$$
x_p (d',n) = x_p (d,n) +1; \; t_p (d',n) = t_p (d,n); \; u_p (d',n) = u_p (d,n),
$$
$$
\al_p (d',n) = \al_p (d,n) + (d+p-1); (\rii \al_p (d-(n+p-1),n) = \al_p (d,n)-d+n);
$$
the same for $\al'_p$}.

b) $\al_{p+1} (d-1,n) \le \al_p (d,n)$, $(\fo) d \ge a^n_p + 1$, $d\in \bz$.

\hspace{0,4 cm} $\al_{p+1} (d,n) \le \al_p (d,n)$, $(\fo) d \ge a^n_p + n+p$, $d\in \bz$.

c) $x_{p+1} (d,n) = x_p (d+1,n+1)$, $t_{p+1} (d,n) = t_p(d+1,n+1)-1$;

\hspace{0,4 cm} $\al_{p} (d+1,n+1) = \al'_{p+1} (d,n) \in \{\al_{p+1}(d,n)-1, \al_{p+1} (d,n)\}$.

{\bf Proof}. b) We'll show first that
$$
\al_{p+1} (d-1,n) \le \al_{p+1} (d,n), \qu (\fo) d \ge a^n_{p+1} + 1, \qu d\in \bz. \leqno (3.3.2)
$$

We use {\it induction} on $x_{p+1} (d-1,n)$, using a).
$$
\left\{\begin{array}{l}
x_{p+1}(d-1,n)=0 \ri \al_{p+1}(d-1,n) = [(d-n+p)/2]_* \\ \noa
x_{p+1}(d-1,n)=0 \ri \al_{p+1}(d,n) = \left\{\begin{array}{l}
[(d-n+p+1)/2]_*, \; d \ne a^n_{p+1} + (n+p) \\ \noa
[(3d-4n+1)/2]_*, \; d = a^n_{p+1} + (n+p).\end{array}\right.
\end{array}\right. \leqno (3.3.3)
$$

So $x_{p+1}(d-1,n)=0 \ri \al_{p+1} (d-1,n) \le \al_{p+1} (d,n)$ (doing some computations for $d=a^n_{p+1}+(n+p))$.

Suppose now that we have proved (3.3.2) for $x_{p+1}(d-1,n) = x\ge 0$. Let $d_1$ so that $x_{p+1} (d_1-1,n) = x+1$. Put $d:= d_1 - (n+p)$. Using a) we deduce that $x_{p+1}(d-1,n) = x$. Then $\al_{p+1}(d_1-1,n) - \al_{p+1}(d_1-1,n) = \al_{p+1}(d+(n+p),n)) - \al_{p+1}
(d-1+(n+p),n) =$ (use a) again) $(\al_{p+1}(d,n) - \al_{p+1}(d-1,n)) + 1 \ge 1 > 0$
(we used the induction hypothesis). Now (3.3.2) is proved.

We'll now prove the first inequality from b) using induction on $x_p (d-1,n)$. We have
$$
x_p(d-1,n) = 0 \ri \al_p (d,n) =
\left\{\begin{array}{l}
[(d-n+p)/2]_*, \; d \ne a^n_{p} + (n+p-1) \\ \noa
[(3d-4n+1)/2]_*, \; d = a^n_{p} + (n+p-1).\end{array}\right. \leqno (3.3.4)
$$

It can be seen that
$$
x_p (d-1,n)=0 \ri [(d-n+p)/2]_* \le \al_p (d,n) \le [(d-n+p+1)/2]_* \leqno (3.3.5)
$$

Because $x_p (d-1,n) =0 \ri x_{p+1} (d-1,n) =0$, from (3.3.3) we now deduce the first inequality from b) if $x_p (d-1,n)=0$.

Suppose now that the inequality holds for $x_p (d-1,n)=x \ge 0$
and let $d_2$ be so that $x_p (d_2-1,n)=x+1$. Put $d:= d_2
-(n+p-1)$. Then  $x_p (d-1,n)=x$ and we have: $\al_p (d_2,n) -
\al_{p+1} (d_2,n) = \al_p (d+(n+p-1),n) - \al_{p+1} (d-1+(n+p),n)
=$ (use a)) $\al_p (d,n) - \al_{p+1} (d-1,n) \ge 0$ (induction
hypothesis). So $\al_p (d_2,n) \ge \al_{p+1} (d_2,n)$. Using now
(3.3.2) we get the first inequality from b).

The second inequality from b) comes from the first one, as before.

\hh

{\bf Lemma 3.3.2}. {\it  (for $F^{p,n}_d (r)$ see (3.2.3))}

a) $F^{p,n}_{d+(n+p-1)} (r+2) = F^{p,n}_d (r) + (d+p-1)$;

\hspace{0,4 cm} $F^{p,n+1}_d (r) = F^{p+1,n}_{d-1} (r) - \di{1\ov 2}$;

b) i) $\al_p (d,n) = \left[F^{p,n}_d (2(x_p(d,n)+1)) - \di{1\ov 2}\right]_*$;

\hspace{0,4 cm} ii) $\al_{p+1} (d-1,n) = \left[F^{p,n+1}_d (2(x_{p+1}(d-1,n)+1)) \right]_*$.

{\bf Proof}. b) Use a), Lemma 3.3.1 a) and  induction on $x_p (d,n)$.

\hh

{\bf Remark}. The previous lemma shows that the functions $\al_p (d,n)$ are related in a "good" way to the genus formula (see (3.2.2)) and  this will be very important in subsection 3.4. Moreover, we already used this property of the functions $\al_p$ when proving the Theorem A (section 2) in [P2].

\subsection{Invertible sheaves ${\cal L} \in \pic (X^n_p)$ whose \\ \protect $(\deg ({\cal L}), P_a ({\cal L}))_{\cal L}$ cover the domain \\ \protect $A^n_p$ $(n\ge 5, p\ge n/3)$}

We recall that the domains $A^n_p$, $p\ge n/3$, $p\in\bz$ from the $(d,g)$-plane were defined in the section 2, (2.11). Here we consider, for $p\ge n/3$, $p\in\bz$, the extended domains $\tilde A^n_p$:
$$
\tilde A^n_p: \al_{p+1} (d-1,n) \le g\le \al_p (d,n), \; d\ge a^n_p +1, \; d,g \in \bz \leqno (3.4.1)
$$
(for the definition of $\al^n_p$ see (2.3)). Obviously
$$
\tilde A^n_p \supset A^n_p, \; (\fo) p\ge n/3, \; (\fo) n\ge 3, \; p,n \in \bz \leqno (3.4.2)
$$
In this section we prove the following

{\bf Proposition 3.4.1.} {\it Let there be $n\ge 5$, $n\in\bz$ and $(d,g) \in \tilde A^n_p$, $p\ge n/3$, $p\in\bz$. Then there is ${\cal D} \in \pic (X^n_p)$ such that $(\deg ({\cal D})$, $p_a ({\cal D})) = (d,g)$.}

We recall that $\deg ({\cal D}):= ({\cal D} \cdot {\cal H}^n_p)$, ${\cal H}^n_p = (p+2; p,1^{s^n_p}) \in \pic (X^n_p)$.

{\bf Proof}. It is divided in 4 steps.

{\it Step 1}. Let there be $d,g \in \bz$, $d\ge a^n_{p+1} +1$ so that
$$
\al'_{p+1} (d-1,n) \le g \le \beta_{p+1} (d,n) - x_{p+1} (d,n),
$$
were $\al'_p$ is defined in (3.3.1) and
$$
\be_p (d,n) = \left\{\begin{array}{ll}
\al'_{p} (d,n), & d\not\equiv a^n_p \; (\hbox{mod}\: (n+p-1)) \\ \noa
\al_{p} (d,n), & d\equiv a^n_p \; (\hbox{mod}\: (n+p-1)) \end{array} \right.
$$
Then there is ${\cal D} \in \pic (X^{n}_p)$ so that $(\deg ({\cal
D}), p_a ({\cal D})) = (d,g)$.

{\it Step 2}. Let there be $d,g \in \bz$, $d\ge a^{n+1}_p$ so that
$$
\al'_{p} (d,n+1) - x_p (d,n+1) \le g \le \al_{p} (d,n+1).
$$
Then there is ${\cal D} \in \pic (X^{n+1}_p)$ so that $(\deg
({\cal D}), p_a ({\cal D})) = (d,g)$.

{\it Step 3}. Let there be $d,g \in \bz$, $d\ge a^{n+1}_p-1 = a^n_{p+1}$ so that
$$
\al'_{p+1} (d,n) - x_{p+1} (d,n) \le g \le \al'_{p+1} (d,n).
$$
Then there is ${\cal D} \in \pic (X^{n}_p)$ so that $(\deg ({\cal
D}), p_a ({\cal D})) = (d,g)$.

{\it Step 4}. The statement of Prop. 3.4.1.

{\bf Proof of Step 1}.  We use induction on $x_{p+1} (d-1,n)$.

If $x_{p+1}(d-1,n)=0$, we have:
$$
\begin{array}{l}
[(d-n+p-1)/2]_* = \al'_{p+1} (d-1,n) \le \be_{p+1} (d,n) - x_{p+1} (d,n) = \\ \noa
= \left\{\begin{array}{ll}
[(d-n+p)/2]_*, & d \ne a^n_{p+1} + (n+p) \\ \noa
[(3d-4n-1)/2]_*, & d = a^n_{p+1} + (n+p). \end{array} \right. \end{array} \leqno (3.4.3)
$$

We consider the following invertible sheaves ${\cal D}_0 \in \pic (X^n_p)$:
$$
{\cal D}_0 = (g+2; g,1^u, 0^{s^n_p-u}), \; 0 \le u \le s^n_p = 3p-n+5. \leqno (3.4.4)
$$

Then $g=p_a ({\cal D}_0)$ and, if $d: = \deg ({\cal D}_0) = ({\cal D}_0 \cdot {\cal H}^n_p)$ it  result that $g = (d-2p+u-4)/2$. Taking now
$$
u = 3p-n+4, \; 3p-n+3, \; 3p-n+2 \ge 0 \leqno (3.4.5)
$$
and using (3.4.3), one obtains invertible sheaves ${\cal D}_0 \in
\pic (X^n_p)$ {\it having the last component zero} and degrees and
arithmetical genera in the range from Step 1, for $x_{p+1}
(d-1,n)=0$. Because
$$
\left\{\begin{array}{l}
\al'_{p+1} (d-1+ (n+p),n) = \al'_{p+1} (d-1,n) + (d+p-1) \\ \noa
(\be_{p+1} - x_{p+1})(d+ (n+p),n) = (\be_{p+1} - x_{p+1}) (d,n) + (d+p-1) \\ \noa
\end{array} \right. \leqno (3.4.6)
$$
(see Lemma 3.3.1 a)) we obtain from (3.4.3), by induction on $x_{p+1} (d-1,n)$ the inequality $\al'_{p+1}(d-1,n) \le \be_{p+1}(d,n) - x_{p+1}(d,n)$, $(\fo) d\ge a^n_{p+1} +1$.

Let there be $\tilde H^n_p = (p+2; p, 1^{s^n_p-1},0) \in \pic (X^n_p)$. It is easy to check that, if ${\cal D} \in \pic (X^n_p)$ {\it has  the last component zero}, then
$$
\left\{\begin{array}{l}
\deg ({\cal D} + \tilde {\cal H}^n_p) = d+(n+p),\; d = \deg ({\cal D}) \\ \noa
p_a ({\cal D} + \tilde {\cal H}^n_p) = g+(d+p-1),\; g = p_a ({\cal D}) \end{array} \right. \leqno (3.4.7)
$$

Now, using (3.4.6) and (3.4.7) it follows that we cover the range from Step 1, adding to the sheaves ${\cal D}_0$ used for $x_{p+1}(d-1,n)=0$, successively, the sheaf $\tilde {\cal H}^n_p$ (at each addition, $x_{p+1} (d-1,n)$ increase by 1, see Lemma 3.3.1 a)).

{\bf Remark 3.4.2.} {\it The invertible sheaves used in order to cover the domain from Step 1 are of the form: ${\cal D} = {\cal D}'_0 + t_2 {\cal H}_2$ (see Lemma 3.2.4), where
${\cal H}_2 = \tilde {\cal H}^n_p = (p+2; p, 1^{s^n_p-1},0)$, $t_2 \in \bz$, $t_2\ge 0$},
$$
{\cal D}'_0 = (a+2; a, 1^u, 0^{s^n_p-u}), \; u \in \{3p-n+2, 3p-n+3, 3p-n+4\},
$$
$$
a\in \{(d'_0-n+p-2)/2, \; (d'_0-n+p-1)/2, \; (d'_0 -n+p)/2\}\cap \bz,
$$
$$
d'_0 = ({\cal D}'_0 \cdot {\cal H}^n_p) \ge a^n_{p+1} + 1, d'_0 \le a^n_{p+1} + (n+p).
$$
(see (3.4.4), (3.4.5) and the construction).

{\bf Proof of Step 2}. This step is the most difficult. We'll give two conctructions, both of them necessary in a complementary way in the subsection 3.5.

{\bf Construction A} (works for $\# \Sigma^n_p \ge 12$, hence $p\ge \di{n\ov 3} + 2$).

1) We'll prove here the existence of ${\cal D} \in \pic (X^{n+1}_p)$ of degree $d$ and arithmetical genus $g = \al_p (d,n+1)$ for any $d\ge a^{n+1}_p$. We recall that
$X^{n+1}_p \su \bp^{n+1}$.

Let there be ${\cal H}^{n+1}_p = (p+2;p,1^{s^{n+1}_p}) \in \pic (X^{n+1}_p)$. Because $\deg (X^{n+1}_p) = n+p$, it can be seen that, for any ${\cal D} \in \pic (X^{n+1}_p)$ we have
$$
\left\{\begin{array}{l}
\deg ({\cal D} + {\cal H}^{n+1}_p) = d+(n+p),\; d = \deg ({\cal D}) = ({\cal D} \cdot {\cal H}^{n+1}_p) \\ \noa
p_a ({\cal D} + {\cal H}^{n+1}_p) = g+(d+p-1),\; g = p_a ({\cal D}). \end{array} \right. \leqno (3.4.8)
$$

Using (3.4.8) and Lemma 3.3.1 a), it follows that we can use the same argument as in the proof of Step 1, doing induction on $x_p (d,n+1)$. So, we check the existence for $x_p (d,n+1)=0$ and then we add successively the class of the hyperplane section of $X^{n+1}_p$, namely ${\cal H}^{n+1}_p$. If $x_p (d,n+1) =0$, then $\al_p (d,n+1)=[(d-n+p-1)/2]_*$. So, we cover this initial range with sheaves
${\cal D}_0 \in \pic (X^{n+1}_p)$ of the form
$$
{\cal D}_0 = (g+2;g,1^u, 0^{s^n_p-u}), \qu g = [(d-n+p-1)/2]_*\; \hbox{with} \leqno (3.4.9)
$$
$$
u = 3p-n+2 = s^n_p -2 \qu \hbox{or} \qu u = 3p-n+3 = s^{n+1}_p -1. \leqno (3.4.10)
$$

Let's {\it remark} that the sheaves used here are of the form
$$
(t+2l; t, l^{s^{n+1}_p-2}, l-1 \; \hbox{or}\; l, l-1),\; \hbox{where} \leqno (3.4.11)
$$
$$
l = x_p(d,n+1) + 1 \in \bz.
$$

2) Now, we consider ${\cal H}_3: = \tilde H^{n+1}_{p-1,1} = (p+1; p-1, 1^{3p-n},0^2,1^2)\in \pic (X^{n+1}_p)$. Then
$$
p_a (\tilde H^{n+1}_{p-1,1}) = p-1, \; \deg (\tilde H^{n+1}_{p-1,1}) = n+p. \leqno (3.4.12)
$$

We start by proving
$$
\left\{\begin{array}{l} \hbox{If}\; g = \al_p (d,n+1) - (x_p
(d,n+1)-x)^2, \; x_p (d,n+1) \ge x, \\ \noa x\in \bz, \; x\ge 0,
\; \hbox{then}\; (\exists) {\cal D} \in \pic (X^{n+1}_p) \;
\hbox{so that} \\ \noa (\deg ({\cal D}), p_a ({\cal D})) = (d,g).
\end{array} \right. \leqno (3.4.13)
$$

We'll use induction on $x_p (d,n+1) \ge x$. During the induction process we'll need the following $\ep ({\cal D})$:
$$
\left\{\begin{array}{l}
\hbox{If}\; {\cal D} = (a; b_0, b_1, \ld, b_{3p-n+4}) \in \pic (X^{n+1}_p) \\ \noa
\hbox{then}\; \ep ({\cal D}):= a- b_0 - b_{3p-n+1} - b_{3p-n+2} \end{array} \right.
\leqno (3.4.14)
$$
(the indices $3p-n+1$ and $3p-n+2$ correspond to the two consecutive zeros from ${\cal H}_3$, so that $\ep ({\cal H}_3)=2$).

If $x_p (d,n+1) =x$, (3.4.13) is just 1). If $x_p (d,n+1) =x$, let's denote by
${\cal D}_x$ an invertible sheaf satisfying (3.4.13), from 1). Then
$\ep ({\cal D}_x) = 0$ $(= 2(x_p(d,n+1)-x))$, by (3.4.11). Actually, we'll prove by induction on $x_p (d,n+1) \ge x$ the following statement (stronger than (3.4.13)):
$$
\left\{\begin{array}{l}
\hbox{If}\; g = \al_p (d,n+1) - (x_p (d,n+1)-x)^2, \; x_p (d,n+1) \ge x \\ \noa
x\in \bz, \; x\ge 0, \; \hbox{then}\; (\exists) {\cal D} \in \pic (X^{n+1}_p) \; \hbox{so that} \\ \noa
(\deg ({\cal D}), p_a ({\cal D})) = (d,g) \; \hbox{and}\; \ep ({\cal D}) = 2(x_p (d,n+1)-x). \end{array} \right. \leqno (3.4.15)
$$

We have just verified (3.4.15) for $x_p (d,n+1) =x$. In order to finish the inductive process, let's suppose that we have constructed invertible sheaves ${\cal D}_{x+t}$ having $(\deg ({\cal D}_{x+t})$, $p_a ({\cal D}_{x+t})) = (d,g)$ in the situation $(t\in \bz,\; t\ge 0)$.
$$
\left\{\begin{array}{l}
x_p (d,n+1) = x + t \\ \noa
g = \al_p (d,n+1) -t^2 = \al_p (d,n+1) - (x_p (d,n+1) -x)^2 \\ \noa
\hbox{so that} \; \ep({\cal D}_{x+t}) = 2t = 2 (x_p (d,n+1) -x)  \end{array} \right. \leqno (3.4.16)
$$

We need ${\cal D}_{x+t+1} \in \pic (X^{n+1}_p)$ to do the same job for the next range in $d$ (namely $(x_p (d',n+1) = x+t+1)$. Or, using the induction hypothesis and Lemma 3.3.1 a) we can see that the sheaves ${\cal D}_{x+t+1}:= {\cal D}_{x+t} + \tilde H^{n+1}_{p-1,1}$ are the good ones. Indeed:
$\deg ({\cal D}_{x+t+1}) = d+(n+p):= d'; x_p (d',n+1) = x+t+1;
p_a ({\cal D}_{x+t+1}) = p_a ({\cal D}_{x+t}) + (p-1) + ({\cal D}_{x+t} \cdot \tilde H^{n+1}_{p-1,1}) -1 = p_a ({\cal D}_{x+t}) + (p-1) + \deg ({\cal D}_{x+t}) - \ep ({\cal D}_{x+t}) -1 =$
(use the induction hypothesis) $\al_p (d,n+1) - (x_p (d,n+1)-x)^2 - 2(x_p (d,n+1)-x) + d+p-2$.

We replace now $d$ by $d'-(n+p)$ and we use Lemma 3.3.1 a). We obtain
$$
\begin{array}{l}
p_a ({\cal D}_{x+t+1}) = \al_p (d'-(n+p),n+1) - (x_p (d'-(n+p),n+1)-x)^2 - \\ \noa
- 2 (x_p (d'-(n+p),n+1)-x) + d'-(n+p)+p-2 = \\ \noa
= \al_p (d',n+1) - d'+n+1 - (x_p (d',n+1)-1-x)^2 - \\ \noa
- 2(x_p (d',n+1)-1-x) + d'-n-2 = \\ \noa
= \al_p (d',n+1) - (x_p (d',n+1)-x)^2 + 2(x_p (d',n+1)-x)-1 - \\ \noa
- 2(x_p (d',n+1)-x)+2-1 = \al_p (d',n+1) - (x_p (d',n+1)-x)^2. \end{array}
$$

Moreover, $\ep ({\cal D}_{x+t+1})  = \ep ({\cal D}_{x+t}) + 2 =
2(x_p (d',n+1)-x)$. Now, the inductive process is finished. So (3.4.15) (hence (3.4.13)) is proved.

Now, let $a$ be {\it arbitrary}, $a\in \bz$, $a\ge 0$ and put $x:=a^2-a\ge 0$. We obtain by (3.4.13) elements from $\pic (X^{n+1}_p)$ of degree $d$ and arithmetical genus
$$
p_a = \al_p (d,n+1) - (x_p (d,n+1) -a^2 + a)^2
$$
for {\it any} $d$ so that $x_p (d,n+1)\ge a^2-a$. Take now $d$ so that
$x_p (d,n+1) = a^2$ $(\ge x = a^2-a)$. One obtains invertible  sheaves of degree $d$ and $p_a = \al_p (d,n+1)-a^2$ for $x_p (d,n+1) = a^2$. Now, adding successively ${\cal H}^{n+1}_p$ and using Lemma 3.3.1 a), we obtain  invertible sheaves of degrees $d$ and arithmetical genera $g$ covering by $(d,g)$ the domain
$$
g = \al_p (d,n+1) = a^2, \qu x_p (d,n+1) \ge a^2. \leqno (a)
$$

3) Now, start again with the invertible sheaves covering the domain ({\it a}) (from 2)). Using the invertible sheaf
$$
{\cal H}_4 = \tilde {\cal H}^{n+1}_{p-1,2} = (p+1; p-1, 1^{3p-n-2}, 0^2, 1^4) \in
\pic (X^{n+1}_p)
$$
(instead of $\tilde H^{n+1}_{p-1,1}$) and using (3.4.11), we obtain, as before, invertible sheaves covering by degrees and arithmetical genera the domains from the $(d,g)$-plane defined by
$$
g =  \al_p (d,n+1) - a^2 - (x_p (d,n+1)-y)^2, x_p (d,n+1)\ge y
$$
(for any $y\ge a^2$, $y\in \bz$).

Take now $b\in\bz$, $b\ge 0$ {\it arbitrary} and put $y = a^2+b^2 -b \ge a^2$. It follows that we obtained invertible sheaves from $\pic (X^{n+1}_p)$ covering by degrees and arithmetical genera the domain
$$
\left\{\begin{array}{l}
g = \al_p (d,n+1) -a^2 - (x_p (d,n+1)-a^2-b^2+b)^2 \\ \noa
x_p (d,n+1) \ge a^2 + b^2 - b. \end{array} \right.
$$

Take now $d$ such that $x_p (d,n+1)=a^2+b^2$ $(\ge y = a^2+b^2-b)$. We obtain invertible sheaves of degrees $d$ and arithmetical genera $g = \al_p (d,n+1) - (a^2+b^2)$. Adding now successively ${\cal H}^{n+1}_{p}$ we get invertible sheaves covering by degrees and arithmetical genera the domain
$$
g = \al_p (d,n+1) - (a^2 + b^2), \; x_p (d,n+1) \ge a^2 + b^2. \leqno (b)
$$

4), 5). Continuity for another two times as for ({\it a}) and ({\it b}), using (3.4.11) (and $s^{n+1}_p \ge 10$, because $p\ge \di{n\ov 3} + 2)$ and the invertible sheaves ${\cal H}_5 =
\tilde {\cal H}^{n+1}_{p-1,3} = (p+1; p-1, 1^{3p-n-4},0^2,1^6)$ and
${\cal H}_6 = \tilde {\cal H}^{n+1}_{p-1,4} = (p+1; p-1, 1^{3p-n-6},0^2,1^8)$,
${\cal H}_5$, ${\cal H}_6 \in \pic (X^{n+1}_p)$ one obtains sheaves from $\pic (X^{n+1}_p)$ covering by degrees and arithmetical genera the domain from the $(d,g)$-plane defined by
$$
g =  \al_p (d,n+1) - (a^2 + b^2 + c^2 + e^2), x_p(d,n+1) \ge a^2 + b^2 + c^2 + e^2 \leqno (e)
$$
where $a,b,c,e \in \bz$, $a,b,c,e \ge 0$ are {\it arbitrary}.

6) The domain ({\it e}) is exactly the domain from Step 2, because any positive integer can be written as a sum of 4 squares of positive integers. Indeed, let there be $g\in \bz$ so that
$$
\al_p (d,n+1) - x_p (d,n+1) \le g \le \al_p (d,n+1)
$$
and
$$
f:= \al_p (d,n+1) - g \ge 0, \qu f\in \bz.
$$
Write $f = a^2 + b^2 + c^2 + e^2$, $a,b,c,e \in \bz$, $a,b,c,e \ge 0$. Certainly $0\le f\le \al_p (d,n+1)$. Hence the domain ({\it e}) from (5) is
$$
\left\{ \begin{array}{l}
g = \al_p (d,n+1) - (a^2 + b^2 + c^2 + e^2) \; (= \al_p (d,n+1)-f) \\ \noa
x_p (d,n+1) \ge a^2 + b^2 + c^2 + e^2 \; (= f = \al_p (d,n+1)-g) \end{array} \right.
$$
$$
0\le f\le x_p (d,n+1).
$$
Replacing $f$ we get $\al_p (d,n+1) - x_p (d,n+1) \le g \le \al_p (d,n+1)$; but this is exactly the domain from Step 2, so it is covered by the $(\deg ({\cal L}), p_a ({\cal L}))$ for various ${\cal L} \in \pic (X^{n+1}_p)$.

{\bf Construction B} (works for $\#\Sigma^n_p \ge 6$, hence $p\ge \di{n\ov 3}$).

From Lemma 3.3.2 b)i) we have
$$
\al_p (d,n+1) = \left[ F^{p,n+1}_d (2(x_p (d,n+1)+1) - {1\ov 2} \right]_* . \leqno (3.4.17)
$$

Due  to this equality we can use the Gruson-Peskine coordinates \\ $(d,r; \te_1, \te_2, \ld, \te_{s^{n+1}_p})$, where $r = a-b_0$, $\te_i = \di{1\ov 2}r - b_i$, $=\overline{1, s^{n+1}_p}$ for ${\cal D} = (a; b_0, b_1, \ld, b_{s^{n+1}_p}) \in \pic (X^{n+1}_p)$ (see (3.2.1)).

Step 2 will be proved as for as for any $(d,g)$ from the domain of the Step 2 we'll prove the existence of a sheaf ${\cal D} \in \pic (X^{n+1}_p)$ with  $\deg {\cal D} =d$, such that, in the coordinates $(d,r; \te_1, \te_2, \ld, \te_{s^{n+1}_p})$, ${\cal D}$ satisfies the conditions (C1) and (C2) from the Prop. 3.2.3 (see [GP2] also), and
$$
F^{p,n+1}_p (r) - {1\ov 2} \sum^{s^{n+1}_p}_{i=1} \te^2_i = g.
\leqno (3.4.18)
$$

Indeed, the left member of (3.4.18) is  just $p_a ({\cal D})$ (see (3.2.2)), hence $(\deg ({\cal D}), p_a ({\cal D})) = (d,g)$. Moreover, it is necessary, performing the transformation of coordinates of ${\cal D}$ from
$$
(d,r; \te_1, \te_2, \ld, \te_{s^{n+1}_p})\; \hbox{to} \; (a; b_0, b_1, \ld, b_{s^{n+1}_p})
$$
that these last coordinates be integers. But (C1) means $b_i \in \bz$, $i = \overline{1, s^{n+1}_p}$ (see (3.2.2)). Now, (C2) means $a\in \bz$, so $b_0 \in \bz$, because $r\in \bz$.

Now we'll prove the existence of such ${\cal D}$ covering the domain from Step 2 by (degree, arithmetical genus). Using (3.4.17), take
$$
r:= 2(x_p(d,n+1)+1). \leqno (3.4.19)
$$

Let there be $b \in \bz$ so that
$$
1\le b \le r =2(x_p(d,n+1)+1)\leqno (3.4.20)
$$
$$
\left\{\begin{array}{l}
\hbox{Write} \; b = \di\sum^4_{i=1} c^2_i, \: c_i \in \bz, \: c_1\ge c_2 \ge c_3 \ge c_4 \ge 0 \\ \noa
\hbox{Then take}\; \te_{s^{n+1}_p - j+1}:= c_j, \: j = \overline{1,4}, \: \te_i = 0,\: i = \overline{s^{n+1}_p-4} \end{array}\right. \leqno (3.4.21)
$$
(recall that $s^{n+1}_p \ge 4$, because $p\ge n/3$).

Given $r$ as in (3.4.19) and $\te_i$ as in (3.4.21), ${\cal D}$ is determined by these numbers and $d$, performing the inverse transformation (3.2.1), with $d = \deg ({\cal D})$.

We have $p_a ({\cal D}) = F^{p,n+1}_d (r) - {1\ov 2} \di\sum^{s^{n+1}_p}_{i=1} \te^2_i =
F^{p,n+1}_d (r) - {1\ov 2}b$ (from 3.2.2) and (3.4.21)). Because $b$ moves in the range (3.4.20), it follows from (3.4.17) that we get the necessary invertible sheaves for Step 2 as far as we check (C1) and (C2). Now, (C1) is clear, because $r\in 2\bz$ and $\te_i \in \bz$. As for  (C2), $g = p_a ({\cal D}) \in \bz$ and $p_a ({\cal D}) = F^{p,n+1}_d (r) - {1\ov 2} \di\sum^{s^{n+1}_p}_{i=1} \te^2_i$; so
$$
2F^{p,n+1}_d (r) - \di\sum^{s^{n+1}_p}_{i=1} \te^2_i \equiv 0 \qu \hbox{(mod 2)}.
$$
Because $l^2 \equiv l$ (mod 2) if $l\in \bz$ and because $r\in 2\bz$, we can see that the left member of the above congruence is congruent (mod 2) with the left member of (C2).

{\bf Remark 3.4.3.} {\it The Construction B is, obviously shorter than the Construction A, but in the subsection 3.5 we'll need both of them}.

{\bf Proof of Step 3.} The existence of the invertible sheaves in the domain from Step 3 in a consequence of the construction from Step 2, using Lemma 3.3.1 c).

Recall that $X^n_p = \va_{[{\cal H}^n_p]} (S^n_p) \su \bp^n$, where $\va_{[{\cal H}^n_p]}$ is the embedding defined by the very ample sheaf ${\cal H}^n_p = (p+2; p, 1^{s^n_p})\in \pic (S^n_p)$ and  $S^n_p = Bl_{\Sigma^n_p} (\bp^2)$ $(\Sigma^n_p \su \bp^2$ a set of $s^n_p+1$ general points). The sheaves from Step 2 were constructed using the surfaces $X^{n+1}_p \su \bp^{n+1}$. If we consider them on the abstract surface $S^{n+1}_p$, there invertible sheaves can be considered on $S^n_p$ also, $S^n_p$ being obtained from
$S^{n+1}_p$ blowing up a new general point $P_{s^n_p}$ (considered in $\bp^2$). Put $\Sigma^n_p:= \Sigma^{n+1}_p \cup \{P_{s^n_p}\}$. The sheaves from $\pic (X^{n+1}_p)$ used in Step 2, considered now in $\pic (S^n_p)$ have the last component equal to zero. Using Lemma 3.3.1c), we obtain invertible sheaves ${\cal D} \in \pic (S^n_p)$ of degree $d$ and $p_a = g$ in the domain
$$
\al'_{p+1} (d-1,n) - x_{p+1} (d-1,n) \le g\le \al'_{p+1} (d-1,n), \; d\ge a^{n+1}_p = a^n_{p+1} + 1.
$$

Now, putting 1 (instead of 0) on the last component of the previous sheaves, the arithmetical genus doesn't change and the degree is translated by 1. We get now exactly the necessary domain for the proof of Step 3.

{\bf Remark 3.4.4.} {\it The invertible sheaves} ${\cal D} \in
\pic (X^n_p)$  {\it used in order to cover the domain from Step 3
comming from the Construction A in Step 2, are of the following
form: ${\cal D} = {\cal D}''_0 + t_2 {\cal H}_2 + t_3 {\cal H}'_3
+ t_4 {\cal H}'_4 + t_5 {\cal H}'_5 + t_6 {\cal H}'_6$, where
${\cal H}_2, {\cal H}'_3, {\cal H}'_4, {\cal H}'_5, {\cal H}'_6$,
are as in Lemma 3.2.4, $t_2, \ld, t_6 \in \bz$, $t_2, \ld, t_6 \ge
0$}
$$
{\cal D}''_0 = (a+2; a, 1^u, 0^{s^{n+1}_p-u}, 1), \; u \in \{3p-n+2, 3p-n+3\}.
$$
$$
a\in \{(d''_0-n+p-1)/2, (d''_0-n+p)/2\} \cap \bz,
$$
$$
d''_0 = ({\cal D}''_0 \cdot {\cal H}^n_p) \ge a^n_{p+1} = [(n-p-1)/2]_*+1,
$$
$$
d''_0 \le a^n_{p+1} +n+p-1.
$$

This follows from (3.4.9), (3.4.10), (3.4.11), the inductive processes, used in proving (a), (b), (c), (e) and the transformation from Step 3 (1 instead of 0 on the last component, giving 1 on the last component of ${\cal D}''_0$ and translating $d_0$ by 1, giving $a$ as in the Remark.)

{\bf Proof of Step 4}. Let's remark at the beginning that, putting together the Steps 1 and 3 and using Lemma 3.3.1 c), it follows that we covered with invertible sheaves from $\pic (X^n_p)$, by degrees and arithmetical genera the domain
$$
\al'_{p+1} (d-1,n) \le g \le \al'_{p+1} (d,n), \qu d\ge a^n_{p+1} +1. \leqno (3.4.22)
$$

We need sheaves from $\pic (X^n_p)$ covering by degrees and arithmetical genera the domain
$$
\al_{p+1} (d-1,n) \le g \le \al_p (d,n), \qu d\ge a^n_p + 1.
$$

We'll do this by induction on $x_p (d-1,n)$, using (3.4.22). If $x_p (d-1,n)=0$, then $x_{p+1} (d-1,n) =0$. Hence $\al_{p+1} (d-1,n) = [(d-n+p)/2]_*$. Using (3.3.5), it follows that we need invertible sheaves which cover by degrees and arithmetical genera the domain
$$
[(d-n+p)/2]_* \le g \le [(d-n+p+1)/2]_*
$$
$$
a^n_p +1 \le d \le a^n_p + (n+p-1).
$$
But then the necessary degrees and genera are realized by
$$
{\cal D}_0 = (g+2; g,1^u, 0^{s^n_p-u}) \in \pic (X^n_p),\; \;
\hbox{with} \leqno (3.4.23)
$$
$$
u = 3p-n+3, \; 3p-n+4 \qu \hbox{or} \qu 3p-n+5. \leqno (3.4.24)
$$

Then, if we succeeded to construct invertible sheaves in the domain $\al_{p+1} (d-1,n) \le g \le \al_p (d,n)$. $x_p(d-1,n)=x$, adding them  the hyperplane class ${\cal H}^n_p = (p+2; p, 1^{s^n_p})$ we obtain, as many times before, using Lemma 3.3.1 a), invertible sheaves from $\pic (X^n_p)$ covering by degrees and arithmetical genera the domain $\al_{p+1} (d,n) \le g \le \al_p (d,n)$ for $x_p (d-1,n) = x+1$. Filling on the left from (3.4.22) (using Lemma 3.3.1 c), last part) we finish the inductive process. Now the Prop. 3.4.1 is proved.

\hh

Before ending this subsection, some remarks, necessary in the next subsection 3.5. The next Remark 3.4.5 is a consequence of the Remarks 3.4.2, 3.4.4 and of the construction given in Step 4.

{\bf Remark 3.4.5.} {\it The invertible sheaves} ${\cal D} \in
\pic (X^n_p)$ {\it used in order to cover by degrees and
arithmetical genera the domain from prop. 3.4.1, using the
Construction A in Step 2, are of one of the following forms:

a) ${\cal D} = {\cal D}'_0 + t_1 {\cal H}_1 + t_2 {\cal H}_2$, where ${\cal H}_1, {\cal H}_2$ are as in the Lemma 3.2.4,} $t_1, t_2 \in \bz$, $t_1,t_2 \ge 0$, ${\cal D}'_0 = (a+2; a, 1^u, 0^{s^n_p-u})$, $s^n_p = 3p-n+5$, $u\in \{3p-n+2, 3p-n+3, 3p-n+4, 3p-n+5\}$, $a\in \{(d'_0-n+p-2)/2$, $(d'_0-n+p-1)/2, (d'_0-n+p)/2, (d'_0-n+p+1)/2\}\cap \bz$,
$d'_0 = ({\cal D}'_0\cdot {\cal H}^n_p) \ge a^n_{p+1} +1$, $d'_0\le a^n_p +n+p$;

b) ${\cal D} = {\cal D}''_0 + t_1 {\cal H}_1 + t_2 {\cal H}_2 + t_3 {\cal H}'_3 +
t_4 {\cal H}'_4 + t_5 {\cal H}'_5 + t_6 {\cal H}'_6$ {\it with ${\cal H}_1, \ld, {\cal H}'_6$ as in the Lemma 3.2.4, $t_1, \ld, t_6 \in \bz$, $t_1, \ld, t_6 \ge 0$, ${\cal D}'_0 = (a+2;a, 1^u, 0^{s^n_p-u},1)$, $u \in \{3p-n+2, 3p-n+3\}$, $a\in \{(d''_0-n+p-1)/2, (d''_0-n+p)/2\}\cap \bz$, $d''_0 = ({\cal D}''_0 \cdot {\cal H}^n_p) \ge a^n_{p+1}$ and} $d''_0 \le a^n_{p+1} +n+p-1$.

From Remark 3.4.5 we deduce the following remark, which we'll use in the next subsection.

{\bf Remark 3.4.6.} {\it If $n\ge 9$ and $\di{n\ov 3} \le p \le
n-4$, $n,p\in \bz$, then the invertible sheaves} ${\cal D} \in
\pic (X^n_p)$ {\it used in order to cover by degrees and
arithmetical genera the domain from Prop. 3.4.1 (using
Construction A in Step 2) are of the following form: ${\cal D} =
{\cal D}_0 + t_1 {\cal H}_1 + t_2 {\cal H}_2 + t_3 {\cal H}'_3 +
t_4 {\cal H}'_4 + t_5 {\cal H}'_5 + t_6 {\cal H}'_6$ with ${\cal
H}_1, {\cal H}_2, {\cal H}'_3, {\cal H}'_4, {\cal H}'_5, {\cal
H}'_6$ as in Lemma 3.2.4 and} ${\cal D}_0 = (a+2;a, 1^u,
0^{s^n_p-u-1},\ep)$, $\ep \in \{0,1\}$, $0\le u\le s^n_p -1 =
3p-n+4$, $a\ge (d_0-n+p-1)/2$, $a\in \bz$, $d_0 \in [a^n_{p+1},
a^n_p +n+p-1]\cap \bz$ $(a^n_p = [(n-p)/2]_*+1)$.

Put $d_0 = d'_0-1$ if ${\cal D}$ is of the category a) in Remark 3.4.5 and $d_0 = d''_0$ if if ${\cal D}$ is of the category b).

{\bf Remark 3.4.7.} {\it Let there be} ${\cal D} = (a; b_0, b_1,
\ld, b_{s^n_p} \in \pic (X^n_p)$ {\it one of the invertible
sheaves used in the proof of Proposition 3.4.1. Let there be $d =
\deg {\cal D} = ({\cal D} \cdot {\cal H}^n_p)$. Then $r\le
2(x_p(d,n)+1)$.}

Indeed, let's denote by $r(d,n):= 2(x_p(d,n)+1)$. If ${\cal D}$  has been used in Step 1, then $r = 2(x_{p+1}(d-1,n)+1)\le r(d,n)$. If ${\cal D}$ has been in Step 3 (with any construction in Step 2), then $r = 2(x_{p+1}(d,n)+1)\le r(d,n)$. If ${\cal D}$  has been used in Step 4, then $r = 2(x_{p+1}(d,n)+1)\le r(d,n)$ or ${\cal D}$ is obtained from some
${\cal D}_1 = (a^1; b^1_0, b^1_1, \ld, b^1_{s_p^n}) \in \pic (X^n_p)$ adding a number, let's say $t$, of ${\cal H}^n_p = (p+2; p, 1^{s^n_p}) \in \pic (X^n_p)$; then $r_1:= a^1 - b^1_0 \le r(d_1,n)$, where $d_1 = \deg ({\cal D}_1)$; now, adding $t$ times ${\cal H}^n_p, r_1$ increases with at most $2t$ in order to become $r$ and $r(d_1,n)$ increases with exactly $2t$ in order to become $r(d,n)$ (due to the form of the three previous $r$); so the inequality $r\le r(d_1,n)$ is transferred to $r\le r(d,n)$ .

\subsection{Curves $C\su X^n_p$ whose $(\deg (C), g(C))_C$ cover
the domain $A^n_p$ $(n\ge 8, n/3 \le p \le n-4)$}

We recall that the domains $A^n_p$ were defined in Section 2, (2.11). In this subsection we prove the following

{\bf Proposition 3.5.1.} {\it Let there be $(d,g)\in A^n_p$, where $n,p\in \bz$, $n\ge 8$, $n/3 \le p \le n-4$, and put $k:= [n/3]_*$. Then:

i) if $d\ge \di{2\ov 3}(3p+n+9)$ there is a (smooth, irreducible) curve $C\su X^n_p$, non-degenerate in $\bp^n$, so that $(\deg (C), g(C)) = (d,g)$;

ii) if $d< \di{2\ov 3}(3p+n+9)$ there is a (smooth, irreducible) curve $C\su \bp^n$, non-degenerate in $\bp^n$, so that $(\deg (C), g(C)) = (d,g)$. If $n\equiv 0$ (mod 3) then $C$ can be found on $X^n_k$ and if $n\equiv 1,2$ (mod 3) then $C$ can be found on $X^n_{k+1}$}.

{\bf Proof}. The proof is a consequence of the analysis which we'll do on the linear systems associated to the invertible sheaves appearing in the proof of Prop. 3.4.1. Precisely, we'll test when these linear systems are nonempty containing (smooth, irreducible) curves, using the smoothing criteria from subsection 3.2. Namely, for the situation (R1) (see (3.1.4)) we'll apply essentially the criterion given by the Corollary 3.2.5 (considering Construction A in Step 2 of the proof of Prop. 3.4.1) and for the situation (R2) (see (3.1.5)) we'll apply the criterion given by the Prop. 3.2.3 (considering Construction B in Step 2 of the proof of Prop. 3.4.1). The situations (R1) and (R2) are complementary and their union cover the hypothesis.

1) {\bf For the situation (R1)}:

{\bf Lemma 3.5.2.} {\it If $(p,n)$ is in the situation (R1), $(d,g) \in A^n_p$ and ${\cal D} \in \pic (X^n_p)$ is one of the invertible sheaves used in the proof of Prop. 3.4.1 (considering Construction A in Step 2) so that $(\deg ({\cal D}), p_a ({\cal D})) = (d,g)$, then we know that ${\cal D} = {\cal D}_0 + t_1 {\cal H}_1 + t_2 {\cal H}_2 + t_3 {\cal H}'_3 + t_4 {\cal H}'_4 + t_5 {\cal H}'_5 + t_6 {\cal H}'_6$ as in Remark 3.4.6.; let's denote by $t:= \di\sum^6_{i=1} t_i$ and suppose $t\ge 3$. Then $[{\cal D}]\ne \emptyset$ and contains a (smooth, irreducible) curve $C$, non-degenerate in $\bp^n$ (recall $\deg ({\cal D}) = ({\cal D} \cdot {\cal H}^n_p))$.}

{\bf Proof}. We use Corollary 3.2.5, i.e. we check that (3.2.10) is verified if $t\ge 3$. Indeed, for $t=3$ (it's enough to consider this case) (3.2.10) becomes
$$
3a \ge 3p-4n+14. \leqno (3.5.1)
$$
Because $a\ge (d_0-n+p-1)/2$, $d_0 \ge (n-p)/2$ (see Remark 3.4.6), minorating $a$ and then $d_0$, it follows that (3.5.1) holds if
$$
13n \ge 9p+62. \leqno (3.5.2)
$$

Now, it is clear that (3.5.2) is satisfied in (R1) $(p\le n-4,\: n\ge 9)$. Because $\deg {\cal H} \ge n+p-1$, $(\fo) {\cal H} \in \{{\cal H}_1, {\cal H}_2, {\cal H}'_3, {\cal H}'_4, {\cal H}'_5, {\cal H}'_6\}$ and $\deg ({\cal D}_0) \ge 0$ (with the definition of the degree as in the statement of the Lemma), the curve $C$ which we just obtained has $\deg (C) \ge 3(n+p-1)$ (because $t= 3$) so $\deg (C) > n+p-1 = \deg (X^n_p)$, and $C\su X^n_p$. So $C$ is non-degenerate in $\bp^n$.

\hh

{\bf Lemma 3.5.3.} {\it If $(p,n), (d,g), {\cal D},t$ are as in the previous lemma and $t=2$, then the same conclusion holds for $[{\cal D}]$}.

{\bf Proof}. We apply again Corollary 3.2.5 for $t=2$ and we use that $a\in \bz$. Doing computations (as before), it follows that the only case when (3.2.10) for $t=2$ is not satisfied is for $p=n-4$ and $a=-1$. We'll consider separately this case.

It's easy to see that the invertible sheaves ${\cal D}$ from our Lemma are of the form
${\cal D} = (b+6; b, \eta_1, \eta_2, \ld, \eta_{s^n_p})$, $\eta_i \in \{0,1,2,3\}$, $(\fo) i$, $b\in \{a+2p-2, a+2p-1, a+2p\}$, $a\in \bz$ as in the Remark 3.4.6. Let $E$ be one of the exceptional divisors lying on $X^n_p$ and corresponding to a $P_i \in \Sigma^n_p$. From the exact sequence
$$
0 \to {\cal D} \to {\cal D} (E) \to {\cal D} (E)|_E \to 0
$$
we deduce that the conclusion of the Lemma holds for $[{\cal D}]$ iff it holds for
$[{\cal D} (E)]$, as far as $b_i = ({\cal D} \cdot {\cal O}_{X^n_p}(E))\ge 1$. Moreover, it is easy to see that all $\eta_i = 3$ only if  $b =a+2p$. So, to conclude the Lemma it remains to prove the conclusion for ${\cal D} \in \{{\cal D}', {\cal D}''\}$ where ${\cal D}' = (a+2p+4; a+2p-2, 3^{s^n_p-1},2)$, ${\cal D}'' = (a+2p+6; a+2p, 3^{s^n_p})$, of course for $p=n-4$ and $a=-1$. We'll study the case ${\cal D} = {\cal D}'$ the other one being similar. So,
$$
{\cal D}' = (2p+3; 2p-3, 3^{2p},2), \qu p=n-4\ge 5.
$$

We specialize the points from $\Sigma^n_p$ on a smooth cubic curve
$\Gamma_0 \su \bp^2$; we denote this specialization by
$\tilde\Sigma^n_p = \{\tilde P_0, \tilde P_1, \ld, \tilde
P_{2p+1}\}$; we suppose, moreover that the points from
$\tilde\Sigma^n_p$ are general on $\Gamma_0$ (see the sketch of
the proof of Prop. 3.2.1), in particular  every 3 of them are not
collinear. Let there by $\tilde S^n_p = Bl_{\tilde\Sigma^n_p}
(\bp^2)$. Then $(p+2;p,1^{2p+1})$ is very ample on $\tilde S^n_p$
(actually, if $s\le 2q+3, q\ge 1, s\ge 0$, $s,q \in \bz$ and
$\Sigma = \{R_0,R_1, \ld, R_s\} \su \Gamma_0$ general on
$\Gamma_0$, $S:=BL_\Sigma (\bp^2)$, then $(q+2; q,1^s)$, is very
ample on $S$, see [Hb2]; for a direct proof see [P4]). Let's
denote by $\tilde X^n_p:= \va_{[p+2;p,1^{2p+1}]} (\tilde S^n_p)
\su \bp^n$. Using Prop. 3.2.1 (i3)$'$ and (ii3)$'$ it results
that, if $\tilde{\cal D}_1:= (2p;2p-3, 2^{2p},1)$, then
$$
\left\{ \begin{array}{l}
h^1 (\tilde{\cal D}_1)=0 \; \hbox{and $[{\cal D}_1] \ne \emptyset$, base point free}, \\ \noa
\hbox{containing a (smooth, irreducible) curve}. \end{array} \right. \leqno (3.5.3)
$$

By semicontinuity (see [CS], Remark 1, p. 324) we deduce that (3.5.3) holds again if we replace $\tilde P_0$ with another point (denoted $\tilde P'_0$) from a small neighborhood of $\tilde P_0$, not belonging to $\Gamma_0$ and not collinear with any others $\tilde P_i$ and $\tilde P_j$. Then, the proper transformation of $\Gamma_0$ in $\tilde X^n_p$ is $\Gamma \in [3;0,1^{2p+1}]$.

Now, we denote by $\de_{ij}$ the quadratic transformation based on $\{\tilde P'_0,
\tilde P_i, \tilde P_j\}$. Performing the successive quadratic transformations $\de_{12}, \de_{34}, \ld, \de_{2p-1,2p}$, the curve $\Gamma_0$ becomes a curve $\dd_0 \su \bp^2$ of degree $p+3$ and having $2p+1$ singular points (with distinct tangents) with multiplicities $p,2,\ld, 2$ respectively; $\tilde X^n_p$ becomes $\bar X^n_p:= \va_{[p+2;p,1^{2p+1}]} (\bar S^n_p) \su \bp^n$, where $\bar S^n_p = Bl_{\bar\Sigma^n_p}
(\bp^2)$, $\bar\Sigma^n_p = \{\bar P_0, \bar P_1, \ld, \bar P_{2p}, \bar P_{2p+1}\}, \bar P_0, \bar P_1, \ld, \bar P_{2p}$ being the singularities of $\dd_0$ with multiplicities $p, \ld, 2$ respectively and $\bar P_{2p+1} \in \dd_0$ a smooth point so that every 3 points in  $\bar X^n_p$ are non-collinear; of course $\tilde X^n_p \cong \bar X^n_p$ $(\cong \tilde S^n_p = \bar S^n_p$); moreover $\tilde {\cal D}_1$ becomes, in the new coordinates, $\bar {\cal D}_1 = (p; p-3, 1^{2p+3}) \in \pic (\bar X^n_p)$. If we consider
$\tilde {\cal D}_1$ and $\bar {\cal D}_1$ on $\tilde S_p = \bar S_p$, then
$\tilde {\cal D}_1 = \bar {\cal D}_1$.

Then, (3.5.3) can be written as
$$
\left\{ \begin{array}{l}
h^1 (\bar{\cal D}_1)=0 \; \hbox{and}\; [\bar{\cal D}_1]\ne \emptyset, \;
\hbox{without base} \\ \noa
\hbox{points and containing a (smooth, irreducible) curve}. \end{array} \right. \leqno (3.5.4)
$$

Let's denote by $\dd \in [p+3;p,2^{2p},1]$ the proper transformation of $\dd_0$ in $\bar X^n_p$. We can see that $\bar{\cal D}_1 = {\cal D}' (-\dd)$, where ${\cal D}'$ is our initial invertible sheaf (considered on $\bar X^n_p$). We obtain then the exact sequence of sheaves on $\bar X^n_p$
$$
0 \to \bar{\cal D}_1 \to {\cal D}'_1 \to {\cal D}'_1|_\dd \to 0.
$$

Since $({\cal D}' \cdot {\cal O}_{\bar X^n_p} (\dd)) = 7\ge 2g(\dd)
= 2 (=2g (\Gamma))$ using (3.5.4) it results that $[\bar{\cal D}']
\ne \emptyset$, base point free containing a (smooth, irreducible)
curve, if the points of $\Sigma^n_p$ are specialized in
$\bar\Sigma^n_p$. So, the same fact remains true on $X^n_p$, by
semicontinuity. \hh

{\bf Lemma 3.5.4.} {\it If $(p,n), (d,g)$, ${\cal D}$ are as in the Lemma 3.5.2 and $d = \deg ({\cal D})\ge \max (2n+1, \di{2\ov 3} (3p+n+9))$, then $[{\cal D}] \ne \emptyset$, without base points and contains a (smooth, irreducible) curve $C$, non-degenerate in $\bp^n$.}

{\bf Proof}. Apply Corollary 3.2.5 for $t=1$.

{\bf  Preliminary Conclusion 3.5.5.} {\it If $(d,g) \in A^n_p$ with $(p,n)$ in situation (R1) and $d\ge max (2n+1,\di{2\ov 3} (3p+n+9))$, then there is a (smooth, irreducible) curve $C \su X^n_p$, non-degenerate in $\bp^n$, with $(\deg (C), g(C)) = (d,g)$}.

\bigskip

{\bf 2) For the situation (R2)}.

{\bf Lemma 3.5.6.} {\it If $(p,n)$ is in the situation (R2), $(d,g) \in A^n_p$ and ${\cal D} \in \pic (X^n_p)$ is one of the invertible sheaves used in the proof of Prop. 3.4.1 (considering Construction B in Step 2) so that $(\deg ({\cal D}), p_a ({\cal D})) = (d,g)$, then ${\cal D}$ satisfied the condition (C3) from the Prop 3.2.3 in the coordinates $(d,r; \te_1, \ld, \te_{s^n_p})$ (here $\deg ({\cal D}) = ({\cal D} \cdot {\cal H}^n_p))$.}

{\bf Proof}. The sheaves used in Step 1 satisfy (C3). The invertible sheaves used in Step 3 come from the sheaves used in Step 2 putting 1 (instead 0) on the last component. So, it's enough to verify that the sheaves ${\cal D}$ used in Step 2 apply to
$$
|\te_1| \le \te_2 \le \ld \le \te_{s^{n+1}_p} < {r\ov 2} \; \hbox{for}\; d = \deg ({\cal D}) \ge 2n+2 \leqno (3.5.5)
$$

The inequalities $|\te_1| \le \te_2 \le \ld \le \te_{s^{n+1}_p}$ come from (3.4.21). From $d\ge 2n+2$ we deduce that $r = 2(x_p(d,n+1)+1) \ge 4$ (cf. (3.4.19)). Then $c_i > {r\ov 2} \ri b \ge c^2_i > \di{r^2 \ov 4}$ (see (3.4.21)); but $b\le r$ (see (3.4.20)) and $r\le
\di{r^2 \ov 4}$, that is a contradiction ! So, $c_i \le \di{r \ov 2}$, hence $\te_i \le \di{r\ov 2}$, $i = \overline{1, s^{n+1}_p}$. Moreover, $c_1 = \di{r \ov 2} \rii r=4$. Then we replace (if necessary)
$(c_1, c_2, c_3, c_4) = (2,0,0,0)$ by $(c_1, c_2, c_3, c_4) = (1,1,1,1)$ and we get
$\te_{s^{n+1}_p} < \di{r\ov 2}$. Now, the invertible sheaves used in Step 4 satisfy (C3), because ${\cal H}^n_p$ satisfies this condition.

\hh

{\bf Lemma 3.5.7.} {\it If $(p,n), (d,g)$ and ${\cal D}$ are as in the previous lemma, then ${\cal D}$ satisfies the condition (C4) from Prop. 3.2.3 in the coordinates
$(d,r; \te_1, \te_2, \ld, \te_{s^{n}_p})$.}

{\bf Proof}. (C4) means $b_0 \ge b_1$ if ${\cal D} = (a; b_0,b_1, \ld, b_{s^n_p})$ is written in the usual coordinates on $\pic (X^n_p)$. We will check this (in i)) for the sheaves used in Step 1 from the proof of Prop. 3.4.1 obtained by adding a (finite) number of $\tilde{\cal H}^n_p$ to the initial family ${\cal D}_0$, for the sheaves used in Step 4 of Prop. 3.4.1 obtained adding a (finite) number of ${\cal H}^n_p$ to the
sheaves used in Step 1 (all these are sheaves like the invertible sheaves from Remark 3.4.5 a)). We'll also check (C4) for the sheaves used in Step 2, Construction B (in ii)). It's clear that from these two verifications (i) and (ii), we get the conclusion of the Lemma.

i) We consider sheaves ${\cal D} = {\cal D}_0 + t_1 {\cal H}_1 + t_2 {\cal H}_2$,
${\cal H}_1 = {\cal H}^n_p$, ${\cal H}_2 = \tilde {\cal H}^n_p$, as in Remark 3.4.5. a).
Put $t:= t_1 + t_2$. If $t=1$, because $d = \deg ({\cal D})\ge 2n+1$, it follows that $d_0:= \deg {\cal D}_0 \ge n-p+1$. But, then ${\cal D} = (a+p+4; a+p, 2^u, 1^{s^n_p-u-1}, \ep)$, $\ep \in \{0,1\}$. We need $a+p\ge 2$. Minorating a as in Remark (3.4.5 a)) and $d_0$ from before, it remains to  have $p\ge 3$, which is true in (R2). If $t\ge 2$, then $b_0 \ge b_1$ (i.e. (C4)) becomes $a+tp \ge t+1$. It's enough to consider the case $t=2$ $(a+tp-t-1)$ is an increasing function of $t$). So, we need $a+2p\ge 3$. Minorating $a$ and $d_0$ from Remark 3.4.5 a), it results that we need $p\ge (n+14)/9$, this last inequality being satisfied in (R2).

ii) Now we are going to study the sheaves used in the proof of Step 2 of Prop. 3.4.1 (Construction B). Their degrees are $d\ge 2n+2$, so $x_p (d,n+1) \ge 1$.

1) If $x_p (d,n+1) \ge 3$, we have $-\te_1 + \di\sum^{s^{n+1}_p}_{i=2} |\te_i| = \di\sum^4_{i=1}c_i$ (see (3.4.21)) $\le 2\sqrt{\di\sum^4_{i=1} c^2_i} \le 2\sqrt{b} \le 2\sqrt{2(x_p(d,n+1)+1)}$ (see 3.4.20). Now, it's enough to prove
$$
\begin{array}{l}
x_p (d,n+1) \ge 3 \ri \\ \noa
\rii d-(n-p+2)(x_p (d,n+1) + 3) \ge 2(x_p (d,n+1)+1)  \end{array}\leqno (3.5.6)
$$
(indeed, then $\sqrt{2(x_p(d,n+1)+1)} \le x_p (d,n+1)+1$, so combining (3.5.6) with the previous inequality, we get exactly (C4) for $r = 2(x_p(d,n+1)+1)$, cf. (3.4.19)). As for (3.5.6), we have
$$
x_p (d,n+3) \ge 3 \ri d\ge (7n+5p+2)/2. \leqno (3.5.7)
$$

a) If $p\ge 4$, then
$$
d-(n-p+2)(x_p (d,n+1)+1)-2(x_p(d,n+1)+1) \ge
$$
$$
\ge d -n(x_p(d,n)+1)\ge d - {n\ov 2(n+p)} (2d+n+3p-2).
$$
This last number is $\ge 0$ iff
$$
d\ge n(n+3p-2)/2p. \leqno (3.5.8)
$$
But in (R2), (3.5.8) follows from (3.5.7). We get then (3.5.6).

b) If $p=3$, then $n\in \{8,9\}$. If $n=8$, (3.5.6) becomes $d\ge 9(x_3 (d,9)+1)$, true for $d\ge 27$. But from (3.5.7) it results that $d\ge 37$, so we are ready. Similarly for $n=9$ ((3.5.6) becomes $d\ge 10(x_3 (d,10)+1)$, true for $d\ge 30$; but (3.5.7) $\ri d\ge 40$.

2) If $x_p (d,n+1)=2$. Then
$$
d\ge (5n+3p+2)/2. \leqno (3.5.9)
$$

In this case (C4) becomes
$$
- \te_1 + \sum^{s^{n+1}_p}_{i=2} \te_i \le d-3 (n-p+2). \leqno (3.5.10)
$$
From (3.5.9) it follows $d-3(n-p+2) \ge (9p-n-10)/2\ge 4$ (since $p\ge n/3$, $n\ge 8$,
and $p=(n+1)/3=3$ if $n=8$). Then in order to obtain (3.5.10) it's enough to choose  $\te_1, \ld , \te_{s^{n+1}_p}$ such that $-\te_1 + \di\sum^{s^{n+1}_p}_{i=2} \te_i \le 4$. By (3.4.20) and (3.4.21) it's enough to express each $b\in \{1,2,\ld, 6\}$ as $b = \di\sum^4_{i=1} c^2_i$ so that $\di\sum^4_{i=1} c_i \le 4$ which is, obviously possible.

3) If $x_p (d,n+1)= 1$ (C4) becomes
$$
-\te_1 + \sum^{s^{n+1}_p}_{i=2} \te_i \le d-2 (n-p+2). \leqno (3.5.11)
$$
But $p\ge 3$ and $d\ge 2n+2$, so $d-2(n-p+2)\ge 4$. By (3.4.20) and (3.4.21) it's enough to express each $b\in \{1,2,3,4\}$ as $b = \di\sum^4_{i=1} c^2_i$ so that $c_i \in \{0,1\}$ for all $i$ (hence $\di\sum^4_{i=1} c_i \le 4$), which is possible.

\hh

{\bf Lemma 3.5.8.} {\it If $(p,n),(d,g), {\cal D}$ are as in Lemma 3.5.6 and $d=\deg ({\cal D}) \ge \max \left(2n+1, \di{2\ov 3} (3p+n+9)\right)$, then ${\cal D}$ satisfies the condition (C5) in the coordinates $(d,r;\te_1, \te_2, \ld, \te_{s^n_p})$}.

{\bf Proof}. From Remark 3.4.7 it follows that it's enough to check (C5) for $r=2(x_p(d,n)+1)$. This condition becomes:
$$
d\ge 2(p-1)(x_p (d,n)+1)+2. \leqno (3.5.12)
$$

If we denote by $E_p (d,n):= d-2(p-1)(x_p (d,n)+1)-2$, then (3.5.12) can be written as
$$
E_p (d,n)\ge 0. \leqno (3.5.13)
$$

We can see that, if (3.5.13) is true for integers $d$ so that $x_p (d,n)=x$, then it is true for any integer $d'$ so that $x_p (d',n)\ge x$. Indeed, if $E_p (d,n)\ge 0$ for all $d$ so that $x_p (d,n)=y$, let $d'$ be so that $x_p (d',n)=y+1$ and $d:= d'-(n+p-1)$. Then  $E_p (d',n) = E_p (d,n) + (n-p-1) \ge E_p (d,n)\ge 0$ (use Lemma 3.3.1 a)). The previous statement follows then by induction on $x_p (d,n)$.

Now, if $x_p (d,n)=2$, (3.5.12) becomes
$$
d\ge 6p-4. \leqno (3.5.14)
$$

Since $x_p (d,n)=2$, we have $d\ge (5n+3p-3)/2$. But then (3.5.14) holds in (R2) (use $3p\le n+5$).

So (3.5.12) is true for any $d$ so that $x_p(d,n)\ge 2$. As for $x_p (d,n)=1$, then (3.5.12) becomes $d\ge 4p-2$. This follows now from $d\ge \di{2\ov 3} (3p+n+9)$ in (R2) (use again $3p\le n+5$).

\hh

Because (C1) and (C2) were verified during the Construction B in the proof of Step 2 of Prop. 3.4.1, and (C3), (C4), (C5) were verified in Lemmas 3.5.6, 3.5.7, 3.5.8, we deduce from  Prop. 3.2.3 the following

{\bf Preliminary Conclusion 3.5.9}. {\it If $(d,g) \in A^n_p$ with $(p,n)$ in the situation (R2) and $d\ge \max \left(2n+1, \di{2\ov 3} (3p+n+9)\right)$, then there is a (smooth, irreducible) curve $C\su X^n_p$, non-degenerate in $\bp^n$, with $(\deg (C), g(C)) = (d,g)$.}

{\bf Lemma 3.5.10.} {\it If $n\ge 8$ and $n/3 \le p\le n-4$, $n,p\in \bz$, $(d,g) \in A^n_p$ and $2n+1 \le d \le \di{2\ov 3} (3p+n+9)$, there is a (smooth, irreducible) curve $C \su X^n_k$ if $n\equiv 0$ (mod 3) and $C \su X^n_{k+1}$ if $n\equiv 1,2$ (mod 3), non-degenerate in $\bp^n$, with $(\deg (C), g(C)) = (d,g)$ (here $k = [n/3]_*$, as usual)}.

{\bf Proof}. We use Prop. 3.2.3. We need curves in the range
$$
\al_{p+1} (d-1,n) \le g\le \al_p (d,n), 2n+1 \le d < {2\ov 3} (3p+n+9).
$$
Then $\al_{p+1} (d-1,n) = [(3d-4n-2)/2]_*$, $\al_p (d,n) = [(3d-4n+1)/2]_*$. Because the set $[\al_{p+1} (d-1,n), \al_p (d,n)]\cap \bz$ does not depend on $p$ we construct the necessary curves on $X^n_k$ if $n=3k$ and on $X^n_{k+1}$ if $n=3k+1, 3k+2$. We remark that $F^{p,n}_d (4) = (3d-4n+2)/2$ (see (3.2.3)).

Precisely, we consider linear systems associated to invertible sheaves ${\cal L} \in \pic (X^n_p)$, $p \in \{k, k+1\}$, $k=[n/3]_*$, which in the Gruson-Peskine coordinates are $(d,4;0^{s^n_p-t}, 1^t)$, with $2n+1 \le d < \di{2\ov 3} (3p+n+9)$, $t\in \{1,2,3,4,5\}$ (see (3.2.1) and (3.2.2)). Hence, $r=4$, $\te_i=0$, $i = \overline{0, s^n_p -t}$,
$\te_{s^n_p-j+1} \in \{0,1\}$, $j = \overline{2,5}$, $\te_{s^n_p}=1$, $s^n_p \in \{5,6,7\}$. We can see immediately that the conditions (C1)-(C5) are satisfied (use $d\ge 2n+1$ and (3.2.2)).

\hh

Now the {\bf proof} of the Prop. 3.5.1 follows immediately from and Lemma 3.5.10.

\hh

\subsection{Conclusion}

Now, the {\bf Proof} of the Theorem C (Section 2) follows, putting together the domains $A^n_p$ from Prop. 3.5.1. Then our Main Theorem follows from (2.15), Theorem A, Theorem B and Theorem C (Section 2).

\hh

\section{Comments and further developments; $D^n_2$ - the expected lacunary domain}

\hspace*{.55cm}1) One can see that our Main Theorem remains true over an albegraically closed field of {\it arbitrary characteristic}. This follows replacing the Bertini Theorem (which we used several times during our proof) - true in characteristic zero - by the Hartshorne's Bertini-type theorem ([Ha2], th\'eor\`eme 5.1) - true in arbitrary characteristic. The verifications are similar to Rathmann's verifications for curves in $\bp^4$ and $\bp^5$ ([Ra]).

2) We can consider other topics from Halphen-Castelnuovo theory (see \S1) also. One example is when ${\cal P}=$ linear normality. In [P5] and [P6] we extended the results of Dolcetti and Pareschi ([DP]) from  $\bp^3$ to higher dimensional projective spaces.

3) Moreover, concerning the study of families of non-degenerate (smooth, irreducible) curves from $\bp^n$, we can consider the Hilbert scheme $H^n_{d,g}$ (see \S1). As Kleppe pointed out ([Kl], \S6) our results can be used in order to stand out "good" components of $H^n_{d,g}$ in a big range on $(d,g,n)$.

4) Finally, few words concerning the domain $D^n_2$ (see (2.10)). We proved that there is no gap for $HC(n)$ in $D^n_1$ (see (2.9)) in our Main Theorem. We recall from \S1 that for $n=3,4,5,6,7$ is known that the $(d,g)$-plane is divided in a lacunary and a non-lacunary domain in each case. We {\it conjecture} that $D^n_1$ and $D^n_2$ from this article represents indeed the good definitions for the {\it non-lacunary} $(D^n_1)$ and the {\it lacunary}  $(D^n_2)$ domains for any $n\ge 7$. Although we do not deal is the present article with $D^n_2$, let's do some hints whose aim is to suggest that $D^n_2$ should be a good definition for the lacunary domain.

If $n\ge 7$ and $k: [n/3]_*$ we consider the  subdomains of $D^n_2$ given by:
$$
\begin{array}{l}
D^{'n}_2: A(d,n) < g \le \pi_0 (d,n), \; d\ge 2n+1 \\ \noa
D^{''n}_2: \left\{\begin{array}{ll} \al_{k+1} (d,n), & \hbox{if} \; n \equiv 1,2 (\hbox{mod 3}) \\ \noa \al_k (d,n), & \hbox{if} \; n\equiv 0 (\hbox{mod 3}) \end{array} \right. < g \le A(d,n), \end{array}
$$
$$
2n+1 \le d < d_1 (n)
$$
(see (1.1), (2.1), (2.6), (2.7)). Then $D^n_2 = D^{'n}_2 \cup D^{''n}_2$. If $(d,g) \in  D^{'n}_2$, inspired from the Harris-Eisenbud Conjecture, true for $d>2^{n+1}$ (see [H], ch. III) and using the Horrowitz results from [Ho], \S 1 and the Example from [Ci], we expect that any such a $(d,g)$ which {\it is not} a gap for $HC(n)$ be the $(deg (C), g(C))$ for some smooth, non-degenerate curve $C \su \bp^n$ lying either on a surface $X^n_p$, $\di{n-5\ov 3} \le p < \di{n\ov 3}$ or on a (possibly singular) scroll. In both cases we get gaps (these $X^n_p$ are obtained blowing up too few points from $\bp^2$, at most 3, and on nrolls there is a formula relating $d$ and $g$, see [Ci], section 2, g). If $(d,g)\in D^{''n}_2$, by the same argument ([H], [Ho]), the necessary curves are expected to lie either on $X^n_p$ $\left(\di{n\ov 2} \le p < \di{n\ov 3}\right)$ or on scrolls. We also expect gaps on $X^n_p$ for $d< a$ functions of degree 3/2 in $n$, whose existence could be proved using an argument similar to the case of cubic surfaces from [GP2]. We do not insist more on $D^n_2$, the classification of gaps from it being a rather delicate problem.

\vspace{0,5 cm}

\hspace{5 cm} Address of author:

\hspace{5 cm} Institute of Mathematics

\hspace{5 cm} Romanian Academy of Sciences

\hspace{5 cm} P.O. Box 1-764

\hspace{5 cm} 014700 Bucharest, Romania

\hspace{5 cm} e-mail: ovidiu\_pasarescu@yahoo.com


\begin{thebibliography}{99}

\bibitem [C] GG.CASTELNUOVO, {\it Sui multipli di una serie lineare di gruppi di punti appartenenti ad una curva algebrica}, Rend. Circ. Matem. Palermo 7(1893).

\bibitem [Ci] CC.CILIBERTO, {\it On the degree and genus of smooth curves in a projective space}, Adv. Math. 81(1990), 198-248.

\bibitem [CS] CC.CILIBERTO, E.SERNESI, {\it Curves on surfaces of degree $2r-\delta$ in $\bp^r$}, Comment. Math. Helv. 64(1989), 300-328.

\bibitem [DP] AA. DOLCETTI, G. PARESCHI, {\it On linearly normal spaces curves}, Math. Z., 198 (1988), 73-82.

\bibitem [Gi] AA.GIMIGLIANO, {\it Our thin knowledge of fat points}, Queen's
Papers in Pure and Appl.Math., 83(1989).

\bibitem [GP1] LL.GRUSON, C.PESKINE, {\it Genre des courbes de l'espace projectif}, Algebraic Geometry, Tromso 1977, LNM 687(1978), 31-59.

\bibitem [GP2] LL.GRUSON, C.PESKINE, {\it Genre des courbes de l'espace projectif (II)}, Ann. Sci. Ec. Norm. Sup., 4-2\`eme s\' erie, 15(1982), 401-418.

\bibitem [Hl] GG.HALPHEN, {\it M\'emoire sur la classification des courbes gauches alg\'ebrique}, J.Ec.Polyt. 52(1882), 1-200.

\bibitem [H] JJ.HARRIS (with colab. of D.EISENBUD), {\it Curves in projective space}, S\' em. Math. Sup. Univ. Montr\'eal, 1982.

\bibitem [Hi] AA.HIRSCHOWITZ, {\it Une conjecture sur la coholomogie des diviseurs sur les surfaces rationelles g\'en\'eriques}, J. Reine Angew. Math., 475(1996), 77-102.

\bibitem [Ha1] RR.HARTSHORNE, {\it Algebraic Geometry}, Berlin-Heidelberg-New York, Springer 1977.

\bibitem [Ha2] RR.HARTSHORNE, {\it Genre des courbes alg\'ebriques dans l'espace projectif (d'apr\'es Gruson-Peskine)}, Sem. Bourbaki, Ast\'erisque 92-93(1982), 301-313.

\bibitem [Hb1] BB.HARBOURNE, {\it Complete linear systems on rational surfaces}, Trans.Amer. Math. Soc., 289(1985), 213-225.

\bibitem [Hb2] BB.HARBOURNE, {\it Very ample divisors on rational surfaces}, Math. Ann. 272(1985), 139-153.

\bibitem [Ho] TT.HOROWITZ, {\it Varieties of law $\Delta$-genus}, Duke Math. J., 50(1983), 667-685.

\bibitem [Kl] JJ.KLEPPE, {\it Concerning the existence nice components in the Hilbert scheme of curves in $\bp^n$ for $n=4$ and 5}, J.reine amgew. Math. 475(1996), 77-102.

\bibitem [Mo] SS.MORI, {\it On degrees and genera of curves on smooth quartic surfaces in $\bp^3$}, Nagoya Math. J, 96(1985), 127-132.

\bibitem [N] MM.NOETHER, {\it Zur Grundlegung der Theorie der algebraischen Raumkurven}, Abhad. der preussischen Ak. der Wiss (1882), 1-120.

\bibitem [P1] OO.P\u AS\u ARESCU, {\it On the existence of algebraic curves in the projective $n$-space}, Arch. Math. 51(1988), 255-265.

\bibitem [P2] OO.P\u AS\u ARESCU, {\it On a theorem of Ciliberto}, An. Univ. Buc. (Mat.-Inf.), L, 2(2001), 133-142.

\bibitem [P3] OO.P\u AS\u ARESCU, {\it On a theorem of Ciliberto and Sernesi}, Alg., Groups and Geom., 20(2003), 301-312.

\bibitem [P4] OO.P\u AS\u ARESCU, {\it Linear systems on some rational surfaces}, Preprint.

\bibitem [P5] OO.P\u AS\u ARESCU, {\it Linearly normal curves in {\bf P}$^4$ and {\bf P}$^5$}, An. St. Univ. Constan\c ta, 9,1 (2001), 73-80.

\bibitem [P6] OO.P\u AS\u ARESCU, {\it Linearly normal curves in {\bf P}$^n$}, Serdica Math. J, 30(2004), 349-362.

\bibitem [Ra] JJ.RATHMANN, {\it The genus of curves in $\bp^4$ and $\bp^5$}, Math. Z. 202(1989), 525-543.

\bibitem [T1] AA.TANNENBAUM, {\it Families of algebraic curves with nodes}, Comp. Math. 41(1980), 107-126.

\bibitem [T2] AA.TANNENBAUM, {\it On the geometric genera of projective curves}, Math. Ann. 240(1979), 213-221.

\bibitem [So] FF.SEVERI, {\it Vorlesungen ueber die algebraische geometrie}, Teubner, Leipzig, 1921.

\end{thebibliography}
\end{document}